\input amstex
\documentstyle{amsppt}
\magnification1200
\tolerance=10000
\overfullrule=0pt
\def\n#1{\Bbb #1}
\def\p{\Bbb C_{\infty}}
\def\r{\Bbb R_{\infty}}

\def\Im{\hbox{Im }}

\def\Hom{\hbox{Hom}}

\def\Ker{\hbox{Ker }}
\def\Lie{\hbox{Lie}}

\def\diag{\hbox{ diag }}

\def\e11{E_{11}}

\def\ga{\goth A}

\def\ve{\varepsilon}
\def\vf{\varphi}

\def\ve{\varepsilon}

\def\vf{\varphi}

\def\de{\delta}

\def\gg{\Gamma}
\def\ga{\gamma}
\def\ka{\kappa}
\def\la{\lambda}

\def\La{\Lambda}

\def\Ga{\Gamma}

\def\be{\beta}

\def\th{\theta}

\def\al{\alpha}

\def\ze{\zeta}
\def\om{\omega}
\def\g{\goth }

\topmatter
\title
Lattice of the dual of an Anderson $t$-motive in terms of a map to a flag variety
\endtitle
\author
A. Grishkov, D. Logachev\footnotemark \footnotetext{E-mails: shuragri{\@}gmail.com; logachev94{\@}gmail.com (corresponding author)\phantom{*******************}}
\endauthor
\thanks Thanks: The authors are grateful to FAPESP, S\~ao Paulo, Brazil for a financial support (process No. 2017/19777-6).
The first author is grateful to SNPq, Brazil, to RFBR, Russia, grant 16-01-00577a (Secs. 1-4), and to Russian Science Foundation,
project 16-11-10002 (Secs. 5-8) for a financial support. The second author is grateful to Greg Anderson, Joseph Bernstein, Vladimir Drinfeld, Urs Hartl and Richard Pink for some important remarks, and to M. Tsfasman for linguistic corrections.
\endthanks
\NoRunningHeads
\address
First author: Departamento de Matem\'atica e estatistica
Universidade de S\~ao Paulo. Rua de Mat\~ao 1010, CEP 05508-090, S\~ao Paulo, Brasil, and Omsk State University n.a. F.M.Dostoevskii. Pr. Mira 55-A, Omsk 644077, Russia.
\medskip
Second author: Departamento de Matem\'atica, Universidade Federal do Amazonas, Manaus, Brasil
\endaddress

\abstract Let $M$ be an uniformizable Anderson $t$-motive of rank $r$, $L$ its lattice and $l_*:=\{l_1,\dots, l_r\}$ its basis. We define a map $\de$ from the set of these bases to a flag variety (the present text gives the definition of $\de$ only for elements of the maximal Schubert cell, and for few other cases). If $l_*$ belongs to the maximal Schubert cell then $\de(l_*)$ is described as a set of matrices parametrized by integer points of a tetrahedron; they are called the Siegel element of $l_*$. We give explicit formulas for a Siegel element for $M'$ --- the dual of $M$. As a by-product, we get another proof of the theorem that the lattice of $M'$ is the dual of the lattice of $M$, independent of the proof obtained by U. Hartl and A.-K. Juschka. Generalizations of this result to non-maximal Schubert cells, other tensor operations (tensor product and Hom) and t-motives having non-trivial endomorphism rings, are subjects of further research. 
\endabstract
\endtopmatter
\document
{\bf 0. Introduction.}
\nopagebreak
\medskip
Paper [GL07] gives a definition of the dual $M'$ of an Anderson t-motive $M$, and proves that the lattice of $M'$ is the dual of the lattice of $M$ --- for the case of $M$ having the nilpotent operator $N=0$. It also contains some explicit formulas for equations defining $M'$.
\medskip
The present paper grew out from an attempt to extend this result to the case of arbitrary $M$ (Theorem 4.5), i.e. of $M$ having the nilpotent operator $N\ne0$. We construct a map $\de$ from the set of bases of lattices of t-motives to a flag variety that we denote by tFV, see 3b.10 for its definition, and 3b.7 for the definition of its Iwahori subgroup. Really, in the present paper the map $\de$ is defined only for the preimage of the maximal Schubert cell of tFV, and for a few more elementary cases where $N=0$ (see Section 3c). Its definition for all cases is a subject of further research. 
\medskip
If $l_*=\{l_1,\dots,l_r\}$ is a basis of a lattice of a t-motive such that $\de(l_*)$ belongs to the maximal Schubert cell of tFV, then $\de(l_*)$ is described in terms of a set of matrices $S_*$ depending on 3 parameters, see (3a.5.4.1). We call it the Siegel object of $l_*$ (it is an analog of a Siegel matrix for the case when $M$ has $N=0$). The set of variation of these parameters is the set of integer points of a 3-dimensional tetrahedron such that each edge contains $\g m$ points, hence the total quantity of matrices $S_*$ is the tetrahedral number $T_{\g m}=\binom{\g m+2}{3}$.
\medskip
These Siegel objects characterize elements in the maximal Schubert cell in a tFV. See for example [B] for the theory of flag varieties. We give in Sections 3a.1 -- 3a.5 a description of these flag varieties and the definition of Siegel objects corresponding to its elements. Also, they can be interpreted as representatives of unilateral cosets in elements of Hecke algebra, see Section 3d.6. The Section 3d.7 contains a "wrong" interpretation in terms of Hecke algebras, and Section 3d.8 contains some research problems. All these sections 3a.1 -- 3d.8 can be read independently on the contents of other parts of the paper. They are of elementary nature, they do not require knowledge of the theory of t-motives. 
\medskip
The main resut of the paper is a description of a Siegel object for $M'$ in terms of a Siegel object of $M$. Namely, the dual Siegel object is a set of matrices $P_*$ (see formula (3.20)) depending on the same parameters. Each of these $P_*$ is a polynomial in some $S_*$. For $P_*$ belonging to a face of the tetrahedron, these polynomials appear in the formula for elements of the inverse of an unitriangular matrix (Remark 3.24B). We do not know any interpretation of other polynomials.
\medskip
The whole proof is mainly combinatorial: we work with quadruple sums involving entries of $S_*$, $P_*$. In fact, the only links relating this combinatorics of Siegel objects with Anderson t-motives, are Propositions 2.2, 2.4 and Lemma 4.1.
\medskip
As a by-product, we get another (independent of the earlier proof of [HJ]) proof of a fact that the lattice of $M'$ is the dual of the lattice of $M$. Since the notations of [HJ] differ from the ones of the present paper, we give in 3.5a the relations between these notations. 
\medskip
{\bf Remark 0.1.} The same Siegel objects appear for a general situation. Let $F_s$ be a "small" subfield of the "big" field $F_b$, $V$ a $n$-dimensional vector space over $F_b$, $N$ a nilpotent operator on $V$ such that $N^{\g m}=0$, $\th\in F_b$ a transcendental element over $F_s$. We denote $R:=F_s[\th\cdot Id_V+N]$. Let $L\subset V$ be a free $R$-module of dimension $r$. It is clear that the Siegel objects describe $L$ in $(V,N)$.
\medskip
Here we describe some possible generalizations of this subject.
\medskip
1. Discrete invariants of $N$ are numbers $k_1,\dots, k_{\g m+1}$ (see (3a.3.2.1)), they come from sizes of Jordan blocks of $N$. For the case $N=0 \ ( \ \iff \g m=1)$ we have $(k_1,k_2)=(r-n,n)$. Anderson t-motives of dimension $n$, rank $r$, having $N=0$ are analogs of abelian varieties of dimension $r$ with multiplication by an imaginary quadratic field, of signature $(r-n,n)$. Hence, the reductive group associated to Anderson t-motives of dimension $n$, rank $r$, having $N=0$ is $GU(r-n,n)$. We can expect that we can associate some object which is a generalization of the reductive group $GU(r-n,n)$, to an Anderson t-motive with invariants $k_1,\dots, k_{\g m+1}$.
\medskip
2. We have $GU(r-n,n)=GL_r\times G_\goth m$ over $\n C$. We can define analogs of Anderson t-motives for other reductive groups, for example $GSp$. What is the analog of the results of the present paper for them? For example, what is the set of Siegel objects?
\medskip
3. We can apply the methods of the paper to tensor products of Anderson t-motives.
\medskip
4. It is known that if $A_1$, $A_2$ are abelian varieties then $A_1\otimes A_2$ is a mixed motive. Is it possible to get an analog of our results for it?
\medskip
5. It is not too difficult to generalize the description of $P_*$ in terms of $S_*$, to the case of non-maximal Schubert cycles. 
\medskip
\medskip
{\bf 1. Definitions.}
Let $q$ be a power of a prime, $\th$ a transcendental element, and $\p$ the completion of an algebraic closure of $\n F_q((1/\theta))$ (topology: $\th^{-n}\to0$). Let $\p[T,\tau]$ be the Anderson ring, i.e. the ring of non-commutative polynomials in two variables $T$, $\tau$ satisfying the following relations (here $a \in \p$):
$$Ta=aT, \ T\tau = \tau T, \ \tau a = a^q \tau \eqno{(1.1)}$$
and $\p[T]$, resp. $\p\{\tau\}$ its subrings of polynomials in $T$, resp. $\tau$,
\medskip
{\bf Definition 1.2.} ([G], 5.4.2, 5.4.18, 5.4.16). An Anderson t-motive\footnotemark \footnotetext{Goss calls these objects abelian t-motives.} $M$ is a left
$\p[T, \tau]$-module which is free and finitely generated both as $\p[T]$- and
$\p\{\tau\}$-module, and such that
$$ \exists \goth m \ \hbox{depending on $M$, such that}\ (T-\theta)^ \goth m M/\tau M=0\eqno{(1.2.1)}$$
\medskip
The dimension of $M$ over $\p\{\tau\}$ (resp. $\p[T]$) is denoted by $n$ (resp.
$r$); this number is called the dimension (resp. rank) of $M$.
\medskip
We shall need an explicit matrix description of t-motives. We denote by $*^t$ the transposition. Let
$e_*=(e_1, ..., e_n)^t$ be the vector column of elements of a basis
of $M$ over $\p\{\tau\}$. There exists a matrix $\goth A\in M_n(\p\{\tau\})$ (depending on $e_*$) such that

$$T e_* = \goth A e_*, \ \ \goth A = \sum_{i=0}^l \goth A_i \tau^i \hbox{ where } \goth A_i
\in M_n(\p). \eqno{(1.3)}$$
Condition (1.2.1) is equivalent to the condition
$$\goth A_0=\theta I_n + N,\eqno{(1.3.1)}$$
where $N$ (depending on $e_*$) is a
nilpotent matrix, and the
condition \{we can choose $ \goth m(M)=1$ \} is equivalent to the condition $N=0$ --- because the $\p$-linear span of $\{e_*\}$ is identified with $M/\tau M$.
\medskip
The Carlitz module $\goth C$ is an Anderson t-motive with $r=n=1$. Let $\{e\}=\{e_1\}$ be the only element of a basis of $M$ over $\p\{\tau\}$.
It is defined uniquely up to the multiplication by an element of $\n F_q^*$. $\goth C$ is given by the equation $Te=\theta e +\tau e$.
We have: $e$ also is the only element of a basis of $\g C$ over $\p[T]$, and the multiplication by $\tau$ is given by $\tau e=(T-\th)e$.
\medskip
The tensor product of Anderson t-motives $M_1$, $M_2$ is defined by $M_1\otimes_{\p[T]}M_2$ where the action of $\tau$ is given by  $\tau(m_1\otimes m_2)=\tau(m_1)\otimes \tau(m_2)$. It is known that $M_1\otimes M_2$ is really a t-motive of rank $r_1r_2$, of dimension $n_1r_2+n_2r_1$. $M_1\otimes M_2$ has $N\ne0$ even if $M_1$, $M_2$ have $N=0$.
\medskip
The $\g m$-th tensor power of $\goth C$ is denoted by $\goth C^{\g m}$. Its rank $r$ is 1 and its dimension is $\g m$.
\medskip
Let $M$ be a t-motive. We fix $\g m$ from (1.2.1). We define the $\g m$-dual of $M$ (denoted by $M^{\prime \g m}$, or simply by $M'$ if $\g m$ is fixed) by the formula

$$M^{\prime \g m}=\Hom_{\p[T]}(M,\goth C^{\g m})$$

where for $\vf\in M^{\prime \g m}$ (i.e., $\vf: M \to \goth C^{\g m}$) the action of $\tau$ on $\vf$ is defined in the standard manner:

$$(\tau(\vf))(m)=\tau(\vf(\tau^{-1}(m)))$$
$M'$ exists not for all t-motives $M$. See [GL07] for a proof that $M'$ exists for all pure $M$ and for large classes of $M$ called standard t-motives ([GL07], Section 11). 
\medskip
{\bf Remark.} Hartl and Juschka in [HJ] consider more general objects (they can be called virtual t-motives). All virtual t-motives have dual. 
\medskip
Let us consider the following exact sequence of left $\p[T,\tau]$-modules $$0\to Z_1 \to Z_2 \to Z_3 \to 0$$ from [G], (5.9.22). $Z_1$ is defined as follows: $$Z_1=\p\{T\}:=\{\sum_{i=0}^\infty a_iT^i\ | \ \underset{i\to\infty}\to{\lim}a_i=0\}$$

For any t-motive $M$ of dimension $n$ and rank $r$ we have an exact sequence of $\n F_q[T]$-modules
$$0\to \Hom_{\p[T,\tau]}(M,Z_1) \overset{\ze}\to{\to} \Hom_{\p[T,\tau]}(M,Z_2)\overset{exp}\to{\to} \Hom_{\p[T,\tau]}(M,Z_3)\eqno{(1.3.1a)}$$
whose terms are denoted by $L(M)$, $\Lie(M)$, $E(M)$ respectively (see [G], 5.9.25, 5.9.17). A fixed basis $e_*$ defines a structure of a $\p$-vector space of dimension $n$ on $E(M)$; this structure depends on $e_*$. Analogically, $e_*$ defines a structure of a $\p$-vector space of dimension $n$ on $\Lie(M)$, but it does not depend on a choice of $e_*$, i.e. $\Lie(M)$ has a canonical structure of $\p$-vector space. 
\medskip
Further, $T$ acts on $\Lie(M)$ as a $\p$-linear operator. We have $N=T-\th I_n$ acts on $\Lie(M)$, it can be identified with the dual of $N$ of (1.3.1), because there is an isomorphism of $\p[T]$-modules $(M/\tau M)^*$ and $\Lie(M)$. Finally, $L(M)$ is a free $\n F_q[T]$-module of dimension $\le r$ (a lattice), and the inclusion $\ze: L(M)\to \Lie(M)$ is a homomorphism of $\n F_q[T]$-modules.
\medskip
Condition \{ exp is surjective \} is equivalent to the condition $\dim_{\n F_q[T]}L(M)=r$ (Anderson; see [G], (5.9.14)). t-motives $M$ satisfying these conditions are called uniformizable. Further on, we consider only uniformizable $M$.

$\ze$ defines an epimorphism $L(M)\underset{\n F_q[T]}\to{\otimes}\p[[T-\th]]\to \Lie(M)$. Its kernel is denoted by $\g q_M$; the exact sequence
$$0\to\g q_M\to L(M)\underset{\n F_q[T]}\to{\otimes}\p[[T-\th]]\to \Lie(M)\to0\eqno{(1.3.2)}$$
is a particular case of the Hodge-Pink structure. It gives us information on the lattice $L(M)\overset{\ze}\to{\hookrightarrow}\Lie(M)$.

\medskip
{\bf Theorem 1.4} (Anderson; see also [P]). Let $M_1$, $M_2$ be two uniformizable t-motives. Then $M_1\otimes M_2$ is also uniformizable, and
$$\goth q_{M_1\otimes M_2}=\goth q_{M_1}\underset{\p[[T-\theta]]} \to{\otimes}\goth q_{M_1}\eqno{(1.4.1)}$$
\medskip
{\bf Remark 1.5.} Uniformizability of $M_1\otimes M_2$ follows from [G], Corollary 5.9.38. A proof of Theorem 1.4 for the case: \{Both $M_1$, $M_2$ have $N=0$\} is given in [GL07], Theorem 6. A proof of Theorem 1.4 can be easily obtained combining the methods of the proof of [GL07], Theorem 6, and of the proof of the present paper.
\medskip
Theorem 1.4 describes the lattice of $M_1\otimes M_2$ in terms of the lattices of $M_1$, $M_2$. Clearly there exists an analog of this theorem for $M'$.
To formulate it, we need to define $\g q'_M$ --- the generalization of the notion of the dual lattice
(see [GL07], Definition 2.3; Section 3) to the case $N\ne0$. Let $L(M)'$ be a free $\n F_q[T]$-module dual to $L(M)$, i.e. such that
the $\n F_q[T]$-pairing $<L(M), L(M)'>$ is perfect. We define $\g q^{\prime \g m}_M=\g q'_M\subset L(M)'\underset{\n F_q[T]}\to{\otimes}\p[[T-\th]]$ as follows:

$$\goth q'_M=\{\ x\in L(M)'\underset{\n F_q[T]}\to{\otimes}\p[[T-\theta]]  $$ $$\hbox{ such that }  \forall y\in \goth q_M  \hbox{ we have } <x,y>\in (T-\theta)^\g m\p[[T-\theta]]\ \} \eqno{(1.6)}$$
\medskip
Roughly speaking, we should have: If $M$ is uniformizable and has dual $M'$ then $M'$ is uniformizable, $L(M')$ is canonically isomorphic to $L(M)'$, and $\g q'_M=\g q_{M'}$. We prove this theorem only for $M$ satisfying Condition 4.2.
\medskip
{\bf 2. Results from the theory of t-motives.}
\medskip
The reader can find in [GL20] a better exposition of the contents of this section. 
\medskip
Let $M$ be an uniformizable t-motive. We fix its $\g m$. We use notations $H^1(M)$, $H_1(M)$ from [G], 5.9.11. Suppose that $M$ has the dual $M'$. \medskip
{\bf Lemma 2.0.} $M'$ is uniformizable.
\medskip
{\bf Proof.} There exists the canonical (up to multiplication by an element of $\n F_q^*$) isomorphism between $H^1(M)$ and $H_1(M')=L(M')$,
see [GL21], (1.9). Further on, we have (see [G], (5.9.14) ):
\medskip
$M$ is uniformizable $\iff h_1(M)=r \ \iff h^1(M)=r \ \iff h_1(M')=r \ \iff M'$ is uniformizable. $\square$
\medskip
{\bf Proposition 2.1.} There is the canonical (up to multiplication by an element of $\n F_q^*$) perfect pairing $L(M)\underset{\n F_q[T]}\to{\otimes}L(M') \to \n F_q[T]$.
\medskip
{\bf Proof.} This is [GL07], Lemmas 5.3.6, 5.3.7. We can also prove existence of this pairing using the composition of the perfect pairing between $H^1(M)$ and $H_1(M)=L(M)$ (see [G], (5.9.35); Goss uses $E=E(M)$ instead of $M$), and the above isomorphism between $H^1(M)$ and $H_1(M')=L(M')$. $\square$
\medskip
{\bf 2.1a.} Let $l_1,\dots, l_r$ be a basis of $L(M)$ over $\n F_q[T]$ and $\vf_1,\dots, \vf_r$ the dual basis of $L(M')$ with respect to this pairing. Further, let $f_1,\dots, f_r$ be a basis of $M$ over $\p[T]$. Since $M'=\Hom_{\p[T]}(M,\goth C^{\g m})$, there exists the dual basis $f'_1,\dots, f'_r$ of $M'$ over $\p[T]$ (we fix a basis element of $\goth C^{\g m}$ which is defined up to multiplication by $\n F_q^*$). We denote these bases by $\hat l$, $\hat \vf$, $\hat f$, $\hat f'$ respectively.
\medskip
Since $L(M)=\Hom_{\p[T,\tau]}(M,Z_1)$, for any $l\in L(M)$ and any $f\in M$ there exists $l(f)\in Z_1=\p\{T\}$. We denote $l(f)$ by $<l,f>_1$ (first pairing). For given bases $\hat l$, $\hat f$, their scattering matrix $\Psi\in M_r(\p\{T\})$ ([A], p. 486) is defined as follows: $\Psi_{ij}:=<l_j,f_i>_1$.
\medskip
Let us consider $\Xi=\sum_{i=0}^\infty
a_iT^i\in\p\{T\}$ of [G], p. 172, line 1; recall that it is the only
element (up to
multiplication by $\n F_q^*$) satisfying
$$\Xi=(T-\theta)\Xi^{(1)}, \hbox{ i.e. }\Xi=(T-\theta)\sum_{i=0}^\infty
a_i^qT^i, \ \ \ \lim_{i\to\infty}a_i=0, \ \ \ |a_0|>|a_i| \ \ \forall
i>0$$ (see [G], p. 171, (*)).
$\Xi$ satisfies the following condition: let $\g C$ be the Carlitz module and $e=f\in \g C$ satisfies $$\tau(f)=(T-\th)f$$ Then $l:=\Xi f$ is a basis element of the lattice of $\g C$, i.e. $\tau(l)=l$.

Equivalently,
$$\Xi=c(1-\th^{-1}T)(1-\th^{-q}T)(1-\th^{-q^2}T)(1-\th^{-q^3}T)...\eqno{(2.1.2)}$$ where $c\in \p$ satisfies $c^{q-1}=-1/\th$. 
\medskip
{\bf Proposition 2.2.} Let $\Psi'$ be the scattering matrix of $M'$ corresponding to the bases $\hat \vf$, $\hat f'$. We have $\Psi'=\Xi^{-\g m}(\Psi^t)^{-1}$.
\medskip
{\bf Proof.} For $\g m=1$ this is [GL07], Lemma 5.4.2. For $\g m>1$ the proof is similar. Let us give it.
We should prove that $$\Psi^t\Psi'=\Xi^{-\g m}I_r\eqno{(2.2.1)}$$ Let $l_i$ be an element of $\hat l$.
We associate it a matrix column $X=X(l_i, \hat f)$, see [GL21], lines above (1.6). $X$ is the $i$-th column of $\Psi$.
Hence, the equality $\Psi^t\Psi'=\Xi^{-\g m}I_r$ is equivalent to the equality $<X(l_i,\hat f), X(\vf_j,\hat f')>=\Xi^{-\g m}\de_{ij}$.

Analogically, for $z\in H^1(M)$ we associate it a matrix line $Y(z)$, see [GL21], lines above (1.6).

Let us consider the explicit formula of the isomorphism $\iota: H_1(M') \to H^1(M)$, see [GL21], proof of Proposition 1.9.
Namely, let $\vf_j\in H_1(M')$. We have $Y(\iota(\vf_j))=\Xi^{-\g m}X(\vf_j,\hat f')^t$ (in [GL21] only the case $\g m=1$ is considered).
The pairing $H_1(M)\underset{\n F_q[T]}\to{\otimes} H^1(M)\to \n F_q[T]$ (see [G], 5.9.35; [GL21], 1.7.1) is defined by [GL21], (1.7.2).
Taking into consideration that $\hat l$, $\hat \vf$ are mutually dual bases with respect to this pairing, we get immediately (2.2.1). $\square$
\medskip
{\bf Remark 2.3.} Both $\Xi$ and the dual bases $\hat \vf$, $\hat f'$ are defined up to multiplication by $\n F_q^*$. Proposition 2.2 should be understood in the form that there is a concordant choice of $\Xi$, $\hat \vf$, $\hat f'$ such that the formula for $\Psi'$ is valid.
\medskip
Now we need a definition of $\th$-shift. Let $\psi=\sum_{i=0}^\infty y_iT^i\in\p\{T\}$. We substitute $T=N+\th$ (here $N$ is an abstract symbol). For some (clearly not for all) $\psi$, as a result of this substitution, we get $\sum_{j=-\ka}^\infty z_{-j}N^j\in\p((T))$ for some $\ka\in \n Z$, $z_{*}\in\p$. In this case, we denote this series $\sum_{j=-\ka}^\infty z_{-j}N^j$ by $\psi_N$. Clearly $\th$-shift is compatible with the multiplication of series. For a scattering matrix $\Psi$ we denote by $\Psi_N$ the result of application of $\th$-shift to all entries of $\Psi$.
\medskip
An elementary calculation shows that $(\Xi^{-1})_N$ exists and that its $\ka$ is equal to 1. Really, (2.1.2) shows that $\Xi^{-1}$ has a simple pole at $T=\th$. 
\medskip
For $x\in \Lie(M)$ and $f\in M$ there exists an element $\partial_x(f)\in\p$, it is canonically defined. See [A] for a definition and [GL20] for explicit formulas in coordinates. 
\medskip
{\bf Proposition 2.4.} For any $l\in L(M)$, $f\in M$ we have: $$(<l,f>_1)_N=\sum_{j=-\g m}^\infty z_{-j}N^j\in\p((T)) \hbox{ exists, has $\ka\le \g m$, and for $i\ge1$ }$$
$$z_i=-\partial_{N^{i-1}(\ze(l))}(f).$$
{\bf Proof.} This is [A], 3.3.2 - 3.3.4. For $\g m=1$ this is [GL07], Lemma 5.6; for $\g m>1$ the proof is similar. Let us give it. We consider two pairings. First, for $l\in L$, $f\in M$ let $<l,f>_1\in Z_1$ be the value of $l(f)$, where $l$ is considered as an element of $\Hom_{\p[T,\tau]}(M,Z_1)$, see (1.3.1a). 
\medskip
Second, for $x\in E(M)$, $f\in M$ let $<x,f>_2\in \p$ be the value of $x(f)$, where $x$ is considered as an element of $E=\Hom_{\p\{\tau\}}(M,\p)$, see [G], proof of 5.6.3, or more explicitly [GL20], Definition 6.2A.
\medskip
According [G], 5.9.25, we have $\Lie(M)=\Hom_{\p[T,\tau]}(M,\n F_q((T^{-1})) \ )$. Hence, for $z\in \Lie(M)$, $f\in M$ an element $z(f)\in \n F_q((T^{-1}))$ is defined. If $z=\ze(l)$ for $l\in L(M)$ then $\ze(l)(f)=<l,f>_1\in Z_1\subset \n F_q((T^{-1}))$. The explicit formula for $<l,f>_1$ is the following:
$$<l,f>_1=\sum_{j=0}^\infty <\exp(T^{-(j+1)}\ze(l)),f>_2T^j\eqno{(2.4.1)}$$
this follows immediately from the results of [G], 5.9; see also [A], p. 486, the first formula of (3.2), or the first formula of the proof of Lemma 3.2.1. Further on, for $z\in \Lie(M)$ we denote $\exp(z)- \alpha(z)$ by $\ve(z)$,
hence
$\sum_{j=0}^\infty <\exp(T^{-(j+1)}\ze(l)),f>_2T^j=\underline{A} +\underline{E}$, where
$$\underline{A}=\sum_{j=0}^\infty <\alpha(T^{-(j+1)}\ze(l)),f>_2T^j; \ \ \ \ \underline{E}=\sum_{j=0}^\infty <\ve(T^{-(j+1)}\ze(l)),f>_2T^j$$
Let us calculate their $\th$-shifts. First, the $\th$-shift of $\underline{E}$ has no terms $N^{j}$ for $j<0$. In fact, we can identify $\Lie(M)$, $E(M)$ with $\p^n$ in such manner that $\exp(z)=\sum_{i=0}^{\infty}C_iz^{(i)}$ where $C_0=I_n$. Hence, we get that
$\ve(z)=\sum_{i=1}^{\infty}C_iz^{(i)}$. This means that for large $j$ the element
$\ve(T^{-(j+1)}\ze(l))$ is small, and hence $\ka(\underline{E})=0$, because finitely many
terms having small $j$ do not contribute to the pole of the $\theta$-shift
of $\underline{E}$ (the reader can prove easily the exact estimations himself, or to look [A],
p. 491).
\medskip
Now, let us consider $\underline{A}_N$. We have
$$T^{-j}=\sum_{i=0}^{\g m-1} (-1)^i \binom{j+i-1}{i}\th^{-(j+i)}N^i$$
hence
$$\underline{A}_N=\sum_{i=0}^{\g m-1} (-1)^i (\sum_{j=0}^\infty \binom{j+i}{i}\th^{-(j+i+1)}T^{j} ) <\al(N^i\ze(l)),f>_2$$
We have
$$(-1)^i \sum_{j=0}^\infty \binom{j+i}{i}\th^{-(j+i+1)}T^{j} = (T-\th)^{-(i+1)},$$
hence $$\underline{A}_N=-\sum_{i=0}^{\g m-1} <\al(N^i\ze(l)),f>_2N^{-(i+1)}$$
This formula implies the proposition.
$\square$
\medskip
{\bf 3a. t-flag varieties.}
\medskip
Let $M$ be a uniformizable t-motive of rank $r$ having $N\ne0$, $L$ its lattice and $l_*:=\{l_1,\dots, l_r\}$ its basis. We shall define a map $\de$ from the set of these bases to a flag variety. We shall call these flag varieties by t-flag varieties (tFV), because of t-motives. 

Since this theory is also valid in characteristic 0, here $\p$ will mean either the previous $\p$ or simply $\n C$. In the latter case, $\th\in \n C$ will mean a transcendental element such that $|\th|>1$, i.e. $\th^{-n}\to0$ if $n\to+\infty$. 
\medskip
We fix a number $n$ and a partition $$n=d_1+...+d_\al\eqno{(3a.1.1)}$$ denoted by $d_*$, where $d_1\ge d_2\ge...\ge d_\al>0$. 
Also, we fix an integer $r$ satisfying $r\ge\al$. 
\medskip
Let $V=\p^n$, $N: V\to V$ a nilpotent operator such that the Jordan form of $N$ consists of 0-Jordan blocks of sizes $d_1, d_2,...,d_\al$ (i.e. there are $\al$ such blocks). Let $\goth m\ge d_1$ be a fixed number, hence $N^\goth m=0$. 
\medskip
Clearly we shall apply these general notations to the case $n, \ r$ are the dimension and rank of a t-motive $M$, $V=Lie(M)$ and $N$ the nilpotent operator of $M$, see (1.3.1). 
\medskip
{\bf Definition 3a.1.2.} The set of $N$-bases $\Cal B=\Cal B(n, d_*, r)$ of type $n, d_*, r$ is the following object. An element of $\Cal B$ is a set of elements $l_1, \dots, l_r\in V$ satisfying the below conditions (3a.1.3.1) -- (3a.1.3.2), up to an  isomorphism $\ga: V\to V$ commuting with $N$, i.e. two sets $l_1, \dots, l_r$ and $l'_1, \dots, l'_r$ define the same element of $\Cal B$ if there exists $\ga$ commuting with $N$ such that $\forall \ i$ we have $\ga(l_i)=l'_i$.
\medskip
Hence, we have: $\dim(\Cal B)=rn$ minus the dimension of the set of $\ga: V\to V$ commuting with $N$, see (3a.3.4) for the explicit formula for this dimension. 
\medskip
Let $T:=\th\cdot Id +N$, where $Id$ is the identity operator. 
\medskip
{\bf 3a.1.3.} $l_1, \dots, l_r$ must satisfy two conditions: 
\medskip
{\bf 3a.1.3.1.} The $\p$-envelope of $T^i(l_j)$ is the whole $V$;
\medskip
{\bf 3a.1.3.2.} $l_j$ (or $T^i(l_j)$ - to check) must be linearly independent over $\n R_\infty$, where $\n R_\infty$ is $\n R$ in characteristic 0, and is $\n F_q((\th^{-1}))$ in characteristic $p$. 
\medskip
{\bf Remark 3a.1.4.} Apparently the condition $r\ge \al$ follows from 3a.1.3 (to check). 
\medskip
{\bf Example 3a.1.5.} Let $\goth m=1$. We have $N=0$, $\al=n$, all $d_i=1$. There exists a map $\de$ from $\Cal B$ to the Grassmannian $Gr(r-n,r)(\p)=Gr(n,r)(\p)$. Its naive definition is the following: let $\g V$ be $\p^r$ with a basis $l_1, \dots, l_r$. Since $l_i\in V$, we have a map $\be: \g V\to V$. (3a.1.3.1) implies that $\be$ is surjective. By definition, $\de(l_1, \dots, l_r)=$ Ker $\be\in Gr(r-n,r)(\p)$ (for the exact definition see Definition 3a.5.4, and 3b.17). 
\medskip
{\bf 3a.2. Origin of the construction: Anderson t-motives and their lattices.} 
\medskip
{\bf Definition 3a.2.1.} A $N$-lattice $L$ of rank $r$ is a subset of $V$ such that there exist elements $l_1, \dots, l_r\in V$ satisfying (3a.1.3) such that $L$ is the set of the sums $$P_1(T)(l_1)+...+P_r(T)(l_r)\eqno{(3a.2.1.1)}$$ where $\forall \ i \ P_i(T)$ is a polynomial in $T$ of some degree $\la_i$: $P_i(T)=\sum_{j=0}^{\la_i}c_{ij}T^j$ where $c_{ij}\in \n F_q$ for finite characteristic, $c_{ij}\in \n Z$ or $\n Q$ (to check!) for characteristic 0. 
\medskip
{\bf Definition 3a.2.2.} Two $N$-lattices $L_1, \ L_2$ are isomorphic if exists an isomorphism $\ga: V\to V$ commuting with $N$ such that $\ga(L_1)=L_2$. 
\medskip
{\bf 3a.2.3.} We see that the set of $N$-lattices up to isomorphism is a quotient of $\Cal B$ by the discrete group $GL_r(\n F_q[\th])$: this is an analog of the fact that a Shimura variety is a quotient of an analog of a Siegel upper half plane $\Cal H_g$ by an analog of $GSp_{2g}(\n Z)$. 
\medskip
{\bf 3a.3. Preliminary definitions; dimension formulas.} 
\medskip
{\bf Definition 3a.3.1.} A partition with zeroes of length $\g r$ of a number $\g n$ is a representation of $\g n$ as a sum
$$\g n = \g d_1+\g d_2+...+\g d_{\g r}$$ where $\g d_i\in \n Z$ and $\g d_1\ge\g d_2\ge...\ge\g d_{\g r}\ge0$, i.e. a partition with zeroes is a partition plus several zeroes at its end. 
\medskip
We extend the partition (3a.1.1) to a partition with zeroes of length $r$:
$$n=d_1+...+d_\al+d_{\al+1}+...+d_r\eqno{(3a.3.2)}$$ where $d_{\al+1}=...=d_r=0$. 

Let $$n=c_1+...+c_{d_1}+c_{d_1+1}+...+c_\goth m$$ be the partition with zeroes of length $\goth m$, conjugate to the partition (3a.3.2) (the definition of the conjugate partition with zeroes of a given length is clear). We have 
\medskip
$\al=c_1\ge c_2\ge...\ge c_{d_1}>0$, $c_{d_1+1}=...=c_\goth m=0$. 
\medskip
Now, we let $c_0:=r$, $c_{\goth m+1}:=0$ and $$\hbox{for }i=1,\dots, \goth m+1 \ \ \ \ k_i:=c_{i-1}-c_i\eqno{(3a.3.2.1)}$$
\medskip
We have $$k_i\ge0, \ \ \ r=\sum_{i=1}^{\goth m+1}k_{i}, \ \ \ \ \ n= k_2+2k_3+3k_4+... + \goth m k_{\goth m+1}\eqno{(3a.3.2.2)}$$ Further, we have $$\dim \Ker N^{i}=c_1+...+c_i,$$ $$\dim N^{\goth m-1}V=k_{\goth m+1},\eqno{(3a.3.3.1)}$$
$$\dim N^{\goth m-2}V=k_\goth m+2k_{\goth m+1},\eqno{(3a.3.3.2)}$$ $$\dim N^{\goth m-3}V=k_{\goth m-1}+2k_\goth m+3k_{\goth m+1}\hbox{ etc., }\eqno{(3a.3.3.3)}$$

{\bf 3a.3.4.} It is easy to check that the dimension of the set of $\ga: V\to V$ commuting with $N$ is $$rn-\sum_{1\le i < j\le \goth m+1}(j-i)k_ik_j\eqno{(3a.3.4.1)}$$ (it does not depend on $r$, but only on $k_2,\dots,k_{\goth m+1}$). See (3d.5.7), (3d.5.7.1). 
\medskip
Does exist a better form of this formula? 
\medskip
{\bf Example 3a.3.5.} For $\goth m=1$ we have: the pair $(k_1,k_2)$ is $(r-n,n)$.
\medskip
Let us also mention formulas (here $i\ge2$)
$$k_{i}:=\dim (\Ker N^{i-1}/ \Ker N^{i-2})- \dim (\Ker N^{i}/ \Ker N^{i-1}) \eqno{(3a.3.6)}$$
$$=\dim (\Im N^{i-2} / \Im N^{i-1}) - \dim (\Im N^{i-1} / \Im N^{i})\eqno{(3a.3.7)}$$
\medskip
{\bf 3a.4. Description of a basis of $N^u(V)$ and of the whole $V$.}
\medskip
Let $l_*=(l_1, \dots, l_r)$ be an element of $\Cal B$. We shall see that $\de(l_*)$ belongs to the maximal Schubert cell of tFV iff the below conditions (3a.4.1), (3a.4.2), (3a.4.4) etc. hold. 
\medskip
Elements of the maximal Schubert cell of tFV are described in terms of some matrices. These matrices are called the Siegel object of $\de(l_*)$, it is a generalization of the Siegel matrix of an element of the maximal Schubert cell of a Grassmannian. Let us describe the Siegel object. 
\medskip
First step: Elements $N^{\goth m-1}l_j$, $j=1,\dots,r$, generate $N^{\goth m-1}V$ as a $\p$-vector space. (3a.3.3.1) shows that the dimension of $N^{\goth m-1}V$ is $k_{\goth m+1}$. 
\medskip
{\bf 3a.4.1.} Condition that $l_*$ belongs to the maximal Schubert cell implies that 
$$\hbox{the last $k_{\goth m+1}$ elements from $l_1,\dots,l_r$ form a $\p$-basis of $N^{\goth m-1}V$.}\eqno{(3a.4.1.1)}$$ 
We denote these elements by $l_{\goth m+1,1}, \dots, l_{\goth m+1,k_{\goth m+1}}$, and their set by $\hat l_{\goth m+1}$. 
\medskip
Further on (second step), elements $N^{\goth m-2}l_j$, $N^{\goth m-1}l_j$, $j=1,\dots,r$, generate $N^{\goth m-2}V$ as a $\p$-vector space. 
\medskip
Elements $N^{\goth m-2}(l_{\goth m+1,i})$, $N^{\goth m-1}(l_{\goth m+1,i})$, $i=1, \dots, k_{\goth m+1}$, are linearly independent over $\p$.  Indeed, let $$\sum_{\al_2=1}^{k_{\g m+1}}c_{\al_2}N^{\g m-2}(l_{\g m+1,\al_2})+ \sum_{\al_1=1}^{k_{\g m+1}}c_{\al_1}N^{\g m-1}(l_{\g m+1,\al_1})=0\eqno{(3a.4.1.2)}$$
be a non-trivial dependence relation. Applying $N$ to (3a.4.1.2) we get $$\sum_{\al_2=1}^{k_{\g m+1}}c_{\al_2}N^{\g m-1}(l_{\g m+1,\al_2})=0.$$ This contradicts (3a.4.1.1). Hence, all $c_{\al_2}$ are 0, and this fact also contradicts (3a.4.1.1).
\medskip
We have (see (3a.3.3.2)) $\dim N^{\goth m-2}V=2k_{\goth m+1}+k_{\goth m}$. Hence, we get:
\medskip
{\bf 3a.4.2.} Condition that $l_*$ belongs to the maximal Schubert cell implies that the last $k_{\goth m}$ elements from the set $l_1,\dots,l_{r-k_{\goth m+1}}$
(we denote them by $l_{\goth m,1}, \dots, l_{\goth m,k_{\goth m}}$ respectively, and their set by $\hat l_{\goth m}$) have the property:
\medskip
{\bf (3a.4.3)} $N^{\goth m-2}(l_{\goth m,i})$, $i=1, \dots, k_{\goth m}$, 
\medskip
$N^{\goth m-2}(l_{\goth m+1,i})$, $i=1, \dots, k_{\goth m+1}$, 
\medskip
$N^{\goth m-1}(l_{\goth m+1,i})$, $i=1, \dots, k_{\goth m+1}$, form a $\p$-basis of $N^{\goth m-2}V$.
\medskip
The third step of the process: elements 
\medskip
$N^{\goth m-3}l_j$, $N^{\goth m-2}l_j$, $N^{\goth m-1}l_j$, $j=1,\dots,r$, 
\medskip
generate $N^{\goth m-3}V$ as a $\p$-vector space. 
\medskip
Elements $N^{\goth m-3}(l_{\goth m,i})$, $N^{\goth m-2}(l_{\goth m,i})$, $i=1, \dots, k_{\goth m}$, 
\medskip
$N^{\goth m-3}(l_{\goth m+1,i})$, $N^{\goth m-2}(l_{\goth m+1,i})$, $N^{\goth m-1}(l_{\goth m+1,i})$, $i=1, \dots, k_{\goth m+1}$, 
\medskip
are linearly independent over $\p$ (the proof is exactly the same as that of the second step). 
\medskip
The quantity of these elements is $2k_\goth m+3k_{\goth m+1}$. We have (see (3a.3.3.3)) $\dim N^{\goth m-3}V=k_{\goth m-1}+2k_\goth m+3k_{\goth m+1}$. 
\medskip
Hence, we get:
\medskip
{\bf 3a.4.4.} Condition that $l_*$ belongs to the maximal Schubert cell implies that the last $k_{\goth m-1}$ elements from the set $l_1,\dots,l_{r-k_{\goth m+1}-k_\goth m}$
(we denote them by $l_{\goth m-1,1}, \dots,\ l_{\goth m-1,k_{\goth m-1}}$ respectively, and their set by $\hat l_{\goth m-1}$) have the property:
\medskip
{\bf (3a.4.5)} $N^{\goth m-3}(l_{\goth m-1,i})$, $i=1, \dots, k_{\goth m-1}$, 
\medskip
$N^{\goth m-3}(l_{\goth m,i})$, $N^{\goth m-2}(l_{\goth m,i})$, $i=1, \dots, k_{\goth m}$, 
\medskip
and $N^{\goth m-3}(l_{\goth m+1,i})$, $N^{\goth m-2}(l_{\goth m+1,i})$, $N^{\goth m-1}(l_{\goth m+1,i})$, $i=1, \dots, k_{\goth m+1}$, 
\medskip
form a $\p$-basis of $N^{\goth m-3}V$.
\medskip
Continuing this process, we represent the ordered set $\{l_1,\dots,l_{r}\}$ as a disjoint ordered union of segments: $$\{l_1,\dots,l_{r}\}=\hat l_1\cup \hat l_2\cup \dots\cup \hat l_{\goth m+1}$$ where the length of $\hat l_{i}$ is $k_i$, i.e. 
$$\hat l_{1}=(l_1, \dots, l_{k_1}), \ \  \hat l_{2}=(l_{k_1+1}, \dots, l_{k_1+k_2}), \ \  \hat l_{3}=(l_{k_1+k_2+1}, \dots, l_{k_1+k_2+k_3})\hbox{ etc.}$$ (we use also a notation 
$$(l_1, \dots, l_{k_1})=(l_{1,1}, \dots, l_{1,k_1}); \ (l_{k_1+1}, \dots, l_{k_1+k_2})=(l_{2,1}, \dots, l_{2,k_2}), \hbox{ etc.)}$$ such that $\forall \ u=0,\dots,\goth m-1$ we have:
\medskip
{\bf (3a.4.6)} A $\p$-basis of $N^uV$ is formed by elements from $N^\al(\hat l_{\be})$, where $\al\in[u,\dots, \goth m-1]$, $\be\in[\al+2,\dots, \goth m+1]$.
\medskip
Explicitly, (3a.4.6) can be described in a form of the diagram 
$$\matrix V&|& \hat l_{2}& \hat l_{3}& \hat l_{4}& \dots&\hat l_{\goth m-1}& \hat l_{\goth m}& \hat l_{\goth m+1}\\ \\
         N(V)  &|&    & N(\hat l_{3})& N(\hat l_{4})& \dots&N(\hat l_{\goth m-1})& N(\hat l_{\goth m})& N(\hat l_{\goth m+1})\\  \\ N^2(V)  &|&
           && N^2(\hat l_{4})& \dots&N^2(\hat l_{\goth m-1})& N^2(\hat l_{\goth m})& N^2(\hat l_{\goth m+1})\\  \\
                \dots                  &|&    \dots       &\dots&\dots&\dots&\dots&\dots&\dots \ \ \  (3a.4.7)\\  \\
N^{\goth m-3}(V)&|&   &&&&N^{\goth m-3}(\hat l_{\goth m-1})& N^{\goth m-3}(\hat l_{\goth m})& N^{\goth m-3}(\hat l_{\goth m+1})\\  \\
N^{\goth m-2}(V)&|&  &&&&& N^{\goth m-2}(\hat l_{\goth m})& N^{\goth m-2}(\hat l_{\goth m+1})\\  \\ N^{\goth m-1}(V)&|& &&&&& & N^{\goth m-1}(\hat l_{\goth m+1})\endmatrix$$

Its meaning is the following. Let us consider the line corresponding to $N^u(V)$. Then elements $N^\al(\hat l_{\be})$ situated on this line and below it, form a basis of $N^u(V)$. 
We see that for any $u$ the set of these $\al, \ \be$ is a triangle. 
\medskip
Some $k_i$ can be 0. The corresponding sets $\hat l_i$ are empty. 
\medskip
\newpage
{\bf 3a.5. Description of the Siegel object of $\de(l_*)$.}
\medskip
Now we can define the Siegel object of $\de(l_*)$. It is a set of Siegel matrices parametrized by points of a tetrahedron. Namely, the sets $$\hat l_1, \ N(\hat l_{2}), \ N^2(\hat l_{3}), \dots, N^{u-1}(\hat l_u),\dots, N^{\goth m-2}(\hat l_{\goth m-1}), \ N^{\goth m-1}(\hat l_{\goth m})\eqno{(3a.5.1)}$$ (they are under the bottom entries of the columns of (3a.4.7) ) 
are linear combinations of all entries of (3a.4.7) which are in the east --- southeast sector from them. 
\medskip
More exactly, parameters of the points of this tetrahedron are integers $u, \ z, \ y$ where $u$ is from (3a.5.1), $y$ is the number of the column of (3a.4.7), and $z$ is the number of the line of (3a.4.7) (we neglect the fact that the numeration starts from 2 for columns, and from 0 for lines). The condition that $N^z\hat l_{y}$ is in the east --- southeast sector from $N^{u-1}(\hat l_u)$, is the following: 
$$u\in [1,\goth m], \ \ z\in [u-1, \goth m-1], \ \  y\in [z+2,\goth m+1]\eqno{(3a.5.2)}$$
It turns out that it is convenient to introduce one more parameter $v$ uniquely defined by $u$:
$$v=u-1\eqno{(3a.5.3)}$$i.e. $N^{u-1}(\hat l_u)=N^{v}(\hat l_u)$. 
\medskip
{\bf Definition 3a.5.4.} The Siegel object of $\de(l_*)$ (case of the maximal Schubert cell of tFV) is the set of Siegel matrices $S_{uvyz}$, where $u, \ v, \ y, \ z$ satisfy (3a.5.2), (3a.5.3), of size $k_u\times k_y$ with entries in $\p$ such that the following holds:
$$N^{v}\hat l_{u}=-\sum_{z=u-1}^{\goth m-1}\sum_{y=z+2}^{\goth m+1}S_{uvyz}\cdot N^z\hat l_{y}\eqno{(3a.5.4.1)}$$
Here $\hat l_{u}, \ \hat l_{v}$ are matrix columns;  if some $k_*$ are 0 then the corresponding $S_{****}$ are empty. The sign minus is not of principle, it comes from applications to the theory of Anderson t-motives. 
\medskip
{\bf 3a.5.5.} Example for $\goth m=3$, $u=1, \ v=0$:
$$\matrix \hat l_{1}=-(S_{1020}\cdot\hat l_{2} &+& S_{1030}\cdot\hat l_{3} &+& S_{1040}\cdot\hat l_{4} \\ \\
&+& S_{1031}\cdot N(\hat l_{3}) &+& S_{1041}\cdot N(\hat l_{4}) \\ \\ &&&+& S_{1042}\cdot N^2(\hat l_{4}) \ )\endmatrix \eqno{(3a.5.5.1)}$$
(terms of a fixed column of this formula correspond to a fixed $y$ and different $z$ of (3a.5.4.1), and terms of a fixed row of this formula correspond to a fixed $z$ and different $y$ of (3a.5.4.1) ).
\medskip
{\bf 3a.5.5.2.} For $\goth m=2$, $u=1, \ v=0$ we have the same diagram (3a.5.5.1) without the rightmost column.
\medskip
\newpage
{\bf 3b. Definition of tFV.}
\medskip
Before giving a definition, we recall a definition of the  "classical" flag varieties. Let $k_*:=\{k_1, \dots,k_{m+1}\}$ and $m$ as above be an ordered sequence of numbers such that $k_i\ge0$, $\sum k_i=r$ (in our situation, $k_i$ are invariants of $N$, see (3a.3.2.1) for their definition). 
\medskip
The "classical" flag variety $X(k_1, \dots,k_{m+1})$ is the set of flags of type $k_1, \dots,k_{m+1}$ in $\p^r$, where a flag of type $k_1, \dots,k_{m+1}$ in $\p^r$ is the set of 
$$0=V_0\subset V_1 \subset V_2 \subset ...  \subset V_m \subset V_{m+1}=\p^r\eqno{(3b.2)}$$where $\dim V_i/V_{i-1}=k_i$. It is well-known that $$\dim X(k_1, \dots,k_{m+1})=\sum_{1\le u<y\le m+1}k_uk_y\eqno{(3b.3)}$$ 

An equivalent definition: let $G=GL_r$, and $P_{k_*}$ be the parabolic subgroup consisting of upper block diagonal matrices, where sizes of diagonal blocks are $k_1, \dots,k_{m+1}$. We have: $$X(k_1, \dots,k_{m+1})=G(\p)/P_{k_*}(\p)\eqno{(3b.4)}$$
Also, let $W$ be an abstract variable and $I_{k_*}\subset G(\p[[W]])$ be the Iwahori subgroup of type $k_*$ defined as follows: let $\pi: G(\p[[W]]) \to G(\p)$ be the natural projection ($W\mapsto0$), then $I_{k_*}:=\pi^{-1}(P_{k_*}(\p))$. We have:
$$X(k_1, \dots,k_{m+1})=G(\p[[W]])/I_{k_*}\eqno{(3b.5)}$$
We can also define the affine flag variety of type $k_*$:
$$Fl_{aff}(k_*):=G(\ \p((W)) \ )/I_{k_*}\eqno{(3b.6)}$$
Typical objects of study are the complete (affine) flag varieties that correspond to $k_*=\{1,\dots, 1\}$. 
\medskip
{\bf Definition 3b.7.} The t-Iwahori subgroup of type $k_*$ (notation: $tI_{k_*}$) is a subgroup of $G(\p[[W]])$ consisting of block matrices such that their diagonal blocks are square blocks of sizes $k_1,\dots,k_{m+1}$, and for any $$1\le u<y\le m+1\eqno{(3b.7a)}$$ all entries of the $(y, \ u)$-th block of this matrix belong to $W^{y-u}\cdot \p[[W]]$. 
\medskip
Example for $m=2$:
$$\left(\matrix &&&|&&&&|\\ &*&&|&&*&&|&*&\\ &&&|&&&&|\\ -&-&-&-&-&-&-&-&-&- \\ &&&|&&&&|\\ W\cdot & \p&[[W]]&|&&*&&|&*&\\ &&&|&&&&|\\ -&-&-&-&-&-&-&-&-&- \\ &&&|&&&&|\\ W^2\cdot & \p&[[W]]&|&W\cdot & \p&[[W]]&|&*&\\ &&&|&&&&|\endmatrix \right)\eqno{(3b.7b)}$$

{\bf Remark 3b.8.} Some $k_i$ can be 0. For example, if $m=4$ and $k_2=0$, then $tI_{k_*}$ is the set of matrices
$$\left(\matrix &&&|&&&&|\\ &*&&|&&*&&|&*&\\ &&&|&&&&|\\ -&-&-&-&-&-&-&-&-&- \\ &&&|&&&&|\\ W^2\cdot & \p&[[W]]&|&&*&&|&*&\\ &&&|&&&&|\\ -&-&-&-&-&-&-&-&-&- \\ &&&|&&&&|\\ W^3\cdot & \p&[[W]]&|&W\cdot & \p&[[W]]&|&*&\\ &&&|&&&&|\endmatrix \right)\eqno{(3b.9)}$$
(sizes of diagonal blocks are $k_1, \ k_3, \ k_4$). 
\medskip
{\bf Definition 3b.10.} The set tFV$(k_*)$ is $G(\p[[W]])/tI_{k_*}$. 
\medskip
Clearly $$\dim_{\p}(tFV(k_*))=\sum_{1\le u<y\le m+1}(y-u)k_uk_y\eqno{(3b.10a)}$$   i.e. it differs from the above formula (3b.2) by coefficients $y-u$. 
\medskip
Representatives of cosets $G(\p[[W]])/tI_{k_*}$ for the maximal Schubert cell are described in Theorem 3d.6.5, see also Example 3d.6.6.
\medskip
We can define the affine analog of tFV$(k_*)$ as follows: tFV$(k_*)_{aff}:=G(\p((W)))/tI_{k_*}$. As we mentioned above, we do not know its relations with t-motives and lattices. 

\medskip
{\bf Remark 3b.11.} A similar, but not exactly the same group, was considered in [R16], Section 2.1. Let $T\subset GL_r$ be the diagonal matrices and $\chi\in X_*(T)^+$. We can identify $\chi$ with a sequence of integers $\al_1\le \al_2 \le ... \le \al_r$. If we fix $k_*$ then a corresponding $\chi$ is defined as follows: $\al_1$ is arbitrary (i.e. $\al_1, \al_2, \dots,  \al_r$ are defined mod $\n Z$), 
\medskip
$\al_1= \al_2 = ... = \al_{k_1}$, 
\medskip
$\al_{k_1+1}= \al_{k_1+2} = ... = \al_{k_1+k_2}=\al_1+1$, 
\medskip
$\al_{k_1+k_2+1}= \al_{k_1+k_2+2} = ... = \al_{k_1+k_2+k_3}=\al_1+2$, 
\medskip
$\dots$
\medskip
$\al_{k_1+k_2+... +k_{m}+1}= \al_{k_1+k_2+... +k_{m}+2} = ... = \al_{r}=\al_1+m$. 
\medskip
{\bf Remark 3b.12.} These $\al_i$ are used in [HJ] in order to describe invariants of $N$, instead of $k_*$. 
\medskip
This $\chi$, and hence $k_*$, defines a vertex in the Bruhat-Tits building of $G(\p[[W]])$ belonging to the standard apartment. Its stabilizer $G_\chi$ is a subgroup of $G(\p[[W]])$ consisting of block matrices such that their diagonal blocks are square blocks of sizes $k_1,\dots,k_{m+1}$, and for any $u,\ y \in [1, m+1]$ all entries of the $(y, \ u)$-th block of this matrix belong to $W^{y-u}\cdot \p[[W]]$. 
\medskip
{\bf 3b.13.} We see that the difference from the Definition 3b.7 is the absence of the condition (3b.7a): $u<y$. Example for $m=2$ (compare with (3b.7b)): 
$$\left(\matrix &&&|&&&&|\\ &\p&[[W]]&|&W^{-1}\cdot & \p&[[W]]&|&W^{-2}\cdot & \p&[[W]]&\\ &&&|&&&&|\\ -&-&-&-&-&-&-&-&-&-&-&- \\ &&&|&&&&|\\ W\cdot & \p&[[W]]&|&&\p&[[W]]&|&W^{-1}\cdot & \p&[[W]]&\\ &&&|&&&&|\\ -&-&-&-&-&-&-&-&-&- &-&-\\ &&&|&&&&|\\ W^2\cdot & \p&[[W]]&|&W\cdot & \p&[[W]]&|&&\p&[[W]]&\\ &&&|&&&&|\endmatrix \right)\eqno{(3b.14)}$$

Really, Richarz considers in [R16], Section 2.1 a more general object than 
\medskip
$\chi$ = a vertex in the Bruhat-Tits building of $G(\p[[W]])$. 
\medskip
Namely, he considers a facet $\g a$. It is clear that its stabilizer $G_\g a$ is a subgroup of the above $G\chi$, but it is never $tI_{k_*}$ from Definition 3b.7. 
\medskip
{\bf Remark 3b.15.} Tensor products. 
\medskip
Let us consider two pairs $(r_1, \chi_1)$ and $(r_2, \chi_2)$, where $\chi_1=(\al_{11}, \al_{12}, \dots,  \al_{1,r_1})$ and $\chi_2=(\al_{21}, \al_{22}, \dots,  \al_{2,r_2})$, and the corresponding tFV denoted by $tFV_1$, $tFV_2$. We denote their $n$'s by $n_1, \ n_2$. We can define their tensor product (denoted by $tFV_1\otimes tFV_2$) as a tFV that corresponds to $(r,\chi)$, where $r=r_1r_2$ and $\chi$ consists of numbers $\al_{1i}+\al_{2j}$, $i=1,\dots, r_1,\ \ \ j=1,\dots,r_2$, ordered in the non-decreasing order. 
\medskip
There is a map $\vf: tFV_1\times tFV_2 \to tFV_1\otimes tFV_2$ defined by the following manner. Let $x_1\in tFV_1, \ x_2\in tFV_2$. We consider the \{$N$-lattices with basis\} $L_1, \ L_2$ corresponding to $x_1, \ x_2$ respectively. We can consider their tensor product $L_1 \otimes L_2$. Recall that it is a $N$-lattice in $\p^n$, where its $n$ is $n_1r_2+n_2r_1$. $L_1 \otimes L_2$ has a distinguished basis. By definition, $\vf(x_1, x_2)$ is an element of $tFV_1\otimes tFV_2$ that corresponds to $L_1 \otimes L_2$ with the above mentioned basis. 
\medskip
The present paper contains the analog of this construction for the functor of duality, see Section 6 for details. Namely, let $(r, \chi)$, where $\chi=(\al_{1}, \al_{2}, \dots,  \al_{r})$ be a pair defining tFV. The dual tFV (denoted by $tFV^*$) is defined by the pair $(r, -\chi)$, where $-\chi:=(-\al_{r}, -\al_{r-1}, \dots,  -\al_{1})$ (if tFV corresponds to $k_*=(k_1, k_2, \dots, k_{m+1})$, then $tFV^*$ corresponds to $k_*=(k_{m+1}, k_m, \dots, k_1$). 
\medskip
We have the duality map $\vf: tFV \to tFV^*$. 
Let $x\in tFV$. If $x$ belongs to the maximal Schubert cell of tFV, then it is defined by a Siegel object, see Section 5. The present paper gives the formula for the Siegel object of $\vf(x)$ in terms of the Siegel object of $x$. This result is used for a proof that for Anderson t-motives the functors of lattice and duality commute. 
\medskip
Analog of this formula for tensor product (i.e. the formula for the Siegel object for $\vf(x_1, x_2)$ in terms of the Siegel objects of $x_1, \ x_2$), as well as analogs of this frmula for $x_1, \ x_2$ that do not belong to the maximal Schubert cell, is a subject of further research. 
\medskip
{\bf Theorem 3b.16.} The set of elements in the maximal Schubert cell of tFV is in 1--1 correspondence with the set of Siegel objects. 
\medskip
{\bf Proof.} Clear. 
\medskip
{\bf 3b.17.} Now we can give a definition of $\de$ for the case of the maximal Schubert cell: $\de(l_*)$ is an element of tFV whose Siegel object is the Siegel object of $l_*$. 
\medskip
{\bf 3c. Case of non-maximal Schubert cells.} 
\medskip
The above construction gives the formulas for $\de$ only for bases $l_*$ such that $\de(l_*)$ belongs to the maximal Schubert cell of tFV. 
\medskip
Let us give some information for the case of other Schubert cells. For simplicity, we shall consider only the case $N=0$. In this case tFV is the Grassmannian $Gr(n,r)$. Let $l_*=\{l_1, \dots, l_r\}$ be a $\n F_q[\th]$-basis of $L$. If $l_1, \dots, l_n$ is a $\p$-basis of $V$, then $\de(l_*)$ belongs to the maximal Schubert cell of $Gr(n,r)$, and $\de(l_*)$ is defined via its Siegel matrix.
\medskip
Let us consider the general case. Let $l_*$ be as above, fixed. There exists a unique sequence of numbers $\al_*:=(\al_1, \dots, \al_n)$ such that $1=\al_1<\al_2< ... < \al_n\le r$ defined as follows: 
\medskip
$\forall \  k\in [1,\dots, n]$ we have: 
\medskip
(1) The dimension of the $\p$-linear envelope of $l_1, \dots, l_{\al_k}$ is $k$;
\medskip
(2) $\forall \  \be\in [\al_k,\dots, \al_{k+1}-1]$ we have: the dimension of the $\p$-linear envelope of $l_1, \dots, l_{\be}$ is $k$ (i.e. $\al_{k+1}$ is the minimal number such that the dimension of the $\p$-linear envelope of $l_1, \dots, l_{\al_{k+1}}$ becomes $k+1$). 
\medskip
{\bf Remark 3c.1.} Since $l_1\ne0$, we have alwais $\al_1=1$. 
\medskip
For Schubert cells in $Gr(n,r)$ we use notations of [B], page 4 ($d, \ n$ of [B] are $n, \ r$ of the present paper). Namely, let $I:=(i_1, \dots, i_n)$ be a multi-index where $1\le i_1<...<i_n\le r$. The set of Schubert cells of $Gr(n,r)$ is in 1 -- 1 with the set of these $I$ (the Schubert cell corresponding to the minimal $I=(1,\dots, n)$ is a point, the Schubert cell corresponding to the maximal $I=(r-n+1,\dots, r)$ is the maximal Schubert cell). 
\medskip
{\bf Proposition 3c.2.} $\de(l_*)$ belongs to the Schubert cell of $Gr(n,r)$ corresponding to $I=(r+1-\al_n, \dots, r+1-\al_1)$.
\medskip
{\bf Remark 3c.3.} Since always $\al_1=1$, for all above $I$ we have $i_n=r$. 
\medskip
{\bf Proof.} To write.
\medskip
{\bf 3c.4. Examples of calculation of $\de$}. 
\medskip
{\bf 1. Case $n=2, \ r=3$.} Let $\al_*= (1,2)$, i.e. the maximal Schubert cell. We have 
$$l_3=s_{11}l_1+s_{12}l_2$$hence $$\de(l_*)=(1:s_{11}:s_{12})\in Gr(2,3)=\n P^2.\eqno{(3c.4.1)}$$Let us consider the change of basis $(l_1, l_2, l_3) \to (l_1, l_3, l_2) $. 
\medskip
Since $$l_2=-\frac{s_{11}}{s_{12}}l_1+\frac{1}{s_{12}}l_3$$ we have $$\de(l_1, l_3, l_2)=(1:-\frac{s_{11}}{s_{12}}:\frac{1}{s_{12}})=(s_{12}:-s_{11}:1)\eqno{(3c.4.2)}$$

Let us apply formulas (3c.4.1), (3c.4.2) for the case $\al_*=(1,3)$, i.e. $l_2=\om l_1$ for some $\om\in \p-\r$. For the basis $(l_1, l_3, l_2) $ we have $s_{11}=\om, s_{12}=0$, i.e. $$\de(l_1, l_3, l_2)=(1:\om:0)$$ and hence (see the linear transformation from (3c.4.2) to (3c.4.1))
$$\de(l_1, l_2, l_3)=(0:-\om:1).$$
This is the only 1-dimensional Schubert cell on $\n P^2$ --- the line at infinity. 
\medskip
{\bf 2. Case $n=2, \ r=4$.} Let $\al_*= (1,2)$, i.e. the maximal Schubert cell. We have 
$$l_3=s_{11}l_1+s_{12}l_2$$
$$l_4=s_{21}l_1+s_{22}l_2$$The matrix of coordinates of two basis vectors of $\de(l_*)$ (considered as a plane in $\p^4$) is $$\left(\matrix 1&0&s_{11}&s_{12}\\0&1&s_{21}&s_{22}\endmatrix \right)\eqno{(3c.4.3)}$$Its $2\times 2$ minors (= Pl\"ucker coordinates) are (here $det=|S|$) $$(1:s_{21}:s_{22}:-s_{11}:-s_{12}:det)\eqno{(3c.4.4)}$$Let us consider the change of basis $(l_1, l_2, l_3, l_4) \to (l_1, l_3, l_2, l_4) $. We have 
$$l_2=-\frac{s_{11}}{s_{12}}l_1+\frac{1}{s_{12}}l_3$$
$$l_4=-\frac{det}{s_{12}}l_1+\frac{s_{22}}{s_{12}}l_3$$hence (3c.4.3), (3c.4.4) for this case are
$$\left(\matrix 1&0&-\frac{s_{11}}{s_{12}}&\frac{1}{s_{12}}\\ \\0&1&-\frac{det}{s_{12}}&\frac{s_{22}}{s_{12}}\endmatrix \right)$$ and $$(-s_{12}:det:-s_{22}:-s_{11}:1:s_{21})\eqno{(3c.4.5)}$$respectively. 
\medskip
Let us apply formulas (3c.4.4), (3c.4.5) for the case $\al_*=(1,3)$, i.e. $l_2=\om l_1$ for some $\om\in \p-\r$. For the basis $(l_1, l_3, l_2,l_4) $ we have $s_{11}=\om, s_{12}=0$, i.e. $$\de(l_1, l_3, l_2,l_4)=(1:s_{21}:s_{22}:-\om:0:\om s_{22})$$ and hence (see the linear transformation from (3c.4.5) to (3c.4.4))
$$\de(l_1, l_2, l_3,l_4)=(0:-\om s_{22}:- s_{22}:-\om:1: s_{21}).$$
This element belongs to the Schubert cell of $Gr(2,4)$ corresponding to $I=(2,4)$.

\medskip
Let us consider the change of basis $(l_1, l_2, l_3, l_4) \to (l_1, l_4, l_2, l_3) $. We have 
$$l_2=-\frac{s_{21}}{s_{22}}l_1+\frac{1}{s_{22}}l_4$$
$$l_3=\frac{det}{s_{22}}l_1+\frac{s_{12}}{s_{22}}l_4$$hence (3c.4.3), (3c.4.4) for this case are
$$\left(\matrix 1&0&-\frac{s_{21}}{s_{22}}&\frac{1}{s_{22}}\\ \\0&1&\frac{det}{s_{22}}&\frac{s_{12}}{s_{22}}\endmatrix \right)$$ and $$(s_{22}:det:s_{12}:s_{21}:-1:-s_{11})\eqno{(3c.4.6)}$$respectively. 
\medskip
Let us apply formulas (3c.4.4), (3c.4.6) for the case $\al_*=(1,4)$, i.e. $l_2=\om_2 l_1$, $l_3=\om_3 l_1$ for some $\om_2,\ \om_3\in \p-\r, \ \om_3/\om_2\not\in\r$. For the basis $(l_1, l_4, l_2, l_3)$ we have: $S=\left(\matrix \om_2&0\\\om_3&0\endmatrix \right)$ and (3c.4.4) is $$(1:\om_3:0:-\om_2:0:0).$$ Hence (see the linear transformation from (3c.4.6) to (3c.4.4))
 we get that for this case $$\de(l_1, l_2, l_3, l_4)=(0:-\om_2:1:0:0:\om_3).$$This element belongs to the Schubert cell of $Gr(2,4)$ corresponding to $I=(1,4)$. 
\medskip
{\bf 3c.5. Further research.} 
\medskip
A basis defines its lattice. Let $\Cal B$ be the set of elements $l_*=(l_1, \dots, l_r)$ such that their $\n F_q[\th]$-linear envelope is a lattice. 
\medskip
$\de$ is a map $\Cal B\to Gr(n,r)(\p)$. 
\medskip
$\de$ is injective (to give a proof), $im \ \de$ is dence in $Gr(n,r)(\p)$. 
\medskip
The group $GL_r(\n F_q[\th])$ acts on $\Cal B$, hence it acts on $im \ \de$. It is an analog of the action of $GL_2(\n Z)$ on the upper half plane in the characteristic 0. 
\medskip
Now we must prove Proposition 3c.2 for all cases, and find its analog for the case $N\ne0$. Further, we want to describe duality and tensor product for non-maximal Schubert cells. 
\medskip
{\bf 3d. Some comments. }
\medskip
{\bf 3d.5.7.}  For any $u$ and $y$ satisfying $1\le u<y\le \goth m+1$ the Siegel object contains $y-u$ matrices of size $k_u\times k_y$, hence the dimension of the set of Siegel objects is
$$\sum_{1\le u<y\le \goth m+1}(y-u)k_uk_y\eqno{(3d.5.7.1)}$$ 
The Theorem 3b.16 implies that this is the dimension of tFV. We have: 
the dimension of the set of elements $l_1, \dots, l_r$ is obviously $rn$. It is easy to check that $\dim GFV$ is $rn$ minus the dimension of the set of $\ga: V\to V$ commuting with $N$, see (3a.3.4.1).
\medskip
{\bf Remark 3d.5.8.} Since always $v=u-1$, in fact, the matrices $S_{uvyz}$ depend on 3 parameters $u,y,z$ satisfying (5.2). Their set is the set of integer points in a tetrahedron. Number $v$ indicates the exponent of $N$ in the left hand side of (5.4.1), by analogy with $z$, which indicates the exponent of $N$ in the right hand side of (5.4.1). This notation is convenient to define a symmetry between $S_*$ and $P_*$, see below.
\medskip
{\bf 3d.6. Description of t-flag varieties in terms of Hecke cosets.}
\medskip
The present section 3d.6 and the next 3d.7, 3d.8 are not logically necessary for the proof. We include them in order to show what there exists another description of generalized flag varieties, and to state some research problems. The reader can skip these sections and continue reading in (3.18). 
\medskip
There are two manners to describe Hecke cosets. Let us consider the simplest model example: $\gg=GL_2(\n Z)$, $\gg_0(p)=\{\left(\matrix a&b\\ c&d \endmatrix \right)\in\gg\ |\ c\equiv 0 \mod p\}$.
\medskip

The first manner: $$\gg=\underset{i\in \n P^1(\n F_p)}\to{\bigcup}\gg_0(p)\cdot\ga_i\hbox{ where for $i\in \n F_p\subset \n P^1(\n F_p)$ we have }\ga_i=\left(\matrix 1&0\\i&1 \endmatrix \right)\eqno{( 3d.6.1)}$$

The second manner: $$\gg\left(\matrix 1&0\\0&p \endmatrix \right)\gg=\underset{i\in \n P^1(\n F_p)}\to{\bigcup}\gg\de_i\hbox{ where for $i\in \n F_p\subset \n P^1(\n F_p)$ we have }\de_i=\left(\matrix 1&i\\0&p \endmatrix \right)\eqno{( 3d.6.2)}$$
\medskip
(in both ( 3d.6.1), ( 3d.6.2)  the set of $i\in \n F_p\subset \n P^1(\n F_p)$ is the affine part = the maximal Schubert cell of $\n P^1(\n F_p)$ ).
\medskip
The second manner is convenient to define multiplication of cosets: if $\de_i, \ \de_j$ are from ( 3d.6.2), then the product of these cosets is $\gg \de_i \de_j$, see any textbook on Hecke algebras.
\medskip
The simplest formula obtained by this manner is: $$T(1,p)^2=T(1,p^2)+pT(p,p)\eqno{( 3d.6.3)}$$
We want to generalize ( 3d.6.1) -- ( 3d.6.3) to other congruence subgroups and to find their relations with GFV. 
\medskip
{\bf Theorem  3d.6.5.} The representatives of the maximal Schubert cell of $G(\p[[W]])/tI_{k_*}$ are block lower unitriangular matrices such that $\forall \ y, \ u$ (as earlier $1\le u<y\le \goth m+1$) all entries of their $(y, \ u)$-th block are polynomials in $W$ of degree $\le y-u$. 
\medskip
{\bf Proof} is straightforward. This is an analog of (3d.6.1). 
\medskip
Example for $\goth m=2$: these representatives are
\medskip
$$\left(\matrix &&&|&&&&|\\ &I_{k_1}&&|&&0&&|&0&\\ &&&|&&&&|\\ -&-&-&-&-&-&-&-&-&- \\ &&&|&&&&|\\ &A_{210}&&|&&I_{k_2}&&|&0&\\ &&&|&&&&|\\ -&-&-&-&-&-&-&-&-&- \\ &&&|&&&&|\\ A_{310}&+&A_{311}W&|&&A_{320}&&|&I_{k_3}&\\ &&&|&&&&|\endmatrix \right)\eqno{( 3d.6.6)}$$
\medskip
where $A_{yu\al}$ is a $(k_y\times k_u)$-matrix with entries in $\p$.
\medskip
Obviously the quantity of these elements is the same as in (3d.5.7.1). 
\medskip
{\bf Question  3d.6.8.} How to describe the simplest isomorphism from $\Cal B$ to $G(\p[[W]])/tI_{k_*}$, at least in terms of $S_{****}$, $A_{***}$ of the maximal Schubert cell? For $\goth m=2$ we identify them by
$$S_{1020}=A_{210}^t, \ \ \ S_{2131}=A_{320}^t$$
but I am not sure that the identification $$S_{1030}=A_{310}^t, \ \ \ S_{1031}=A_{311}^t$$
is the best one. (To formulate better).
\medskip
{\bf Question  3d.6.9.} How to interpret matrices $P_{****}$ from (3.20) in terms of $G(\p[[W]])/tI_{k_*}$?
\medskip
{\bf Question  3d.6.10.} The set of representatives of cosets (3d.6.6) is a (non-abelian) group. If $\goth m>1$ then it does not coincide with the natural group structure on the set of $S_{****}$ defined by matrixwise summation of $S_{****}$.

What group structure is more natural; does exist its interpretation in terms of t-motives?

\medskip
{\bf  3d.7. Analog of the second manner.} 
\medskip
An analog of $\left(\matrix 1&0\\0&p\endmatrix \right)$ is the following.
\medskip
{\bf Notation 3d.7.1.} Let $k_*=(k_1,\dots, k_{\goth m+1})$ be given, $\sum_i k_i=r$. We denote the matrix 
$$\diag(1,\dots,1,W,\dots, W, \dots \dots , W^\goth m, \dots, W^\goth m),$$where $W^i$ appears $k_{i+1}$ times, by $\diag(k_*)$. 
\medskip

{\bf 3d.7.2.} Let us consider the Hecke double coset
$\Ga \diag(k_*)\Ga$. It is a union of right cosets. But the set of these right cosets is not $GFV(n, d_*, r)$ if $\goth m\ge2$ (to check carefully). Really, let us consider the simplest counterexample $\goth m=2$, $k_i=1$. Elements
$$\left(\matrix 1&a_{120}&a_{130}+a_{131}W\\ 0&W&a_{230}+a_{231}W\\ 0&0&W^2 \endmatrix \right)\eqno{(3d.7.3)}$$
(where $a_{***}\in \p$) are representatives of cosets (to check).
\medskip
{\bf 3d.8. Some questions.}
\medskip
This section contains some questions related to tFV. 
\medskip
{\bf 3d.8.1. } We have two cosets $\gg_0(p)\cdot\ga_i$ and $\gg_0(p)\cdot\ga_j$; how to get their product, i.e. a coset of $\gg_0(p^2)$? 
\medskip
We call it the Hecke product construction. Clearly we need it not for the above model example, but for the case of right cosets of $\gg_0(k_1,\dots,k_{\goth m+1})$. First, we define the sum of $k_*=(k_1,\dots,k_{\goth m+1})$ and of $k'_*=(k'_1,\dots,k'_{\goth m+1})$ (both $\sum_ik_i=\sum_ik'_i=r$) as the product of the corresponding diagonal matrices: $$k_*+k'_*=k''_*:=(k''_1,\dots,k''_{\goth m+1})\iff \diag(k_*)\cdot\diag(k'_*)=\diag(k''_*)\eqno{(3d.8.2)}$$

We must have a map (a sum) from $[\gg_0(k_*)\backslash \gg] \times [\gg_0(k'_*)\backslash \gg]$ to $\gg_0(k''_*)\backslash \gg$ (maybe defined not at all elements of $[\gg_0(k_*)\backslash \gg] \times [\gg_0(k'_*)\backslash \gg]$).
\medskip
{\bf Question 3d.8.3.} What is its definition? What is its interpretation in terms of first, $N$-lattices, second, t-motives? 
\medskip
{\bf 3. Matrices $P_{****}$ and formulas for duality.}
\medskip
Before starting this section, we indicate that the objects $F_b$, $F_s$ from Remark 0.1 for our case are $F_b=\p$, $F_s=\n F_q$. Further, for the reader's convenience, we indicate the relations between the objects and notations of [HJ] and the present paper. First, instead of $$L=[\Hom _{\p[T]}(M,Z_1)]^\tau$$ Hartl and Juschka consider the dual $\n F_q[T]$-module $\La:=[M\otimes_{\p[T]}Z_1]^\tau$. They consider virtual Anderson t-motives; the t-motives considered in the present paper are called effective in [HJ]. Further, for an effective $M$ the lattice $\g q$ of [HJ] is the dual of $\g q$ of the present paper. Particularly, it contains $\La\otimes_{\n F_q[T]}\p[[N]]$ while in the present paper $\g q $ is contained in $ L(M)\otimes_{\n F_q[T]}\p[[N]]$. Finally, $N$ of [HJ] is $T^{-1}-\th^{-1}$ and not $T-\th$ like in the present paper. 
\medskip
Let $k_1, \dots, k_{\g m+1}$ be from (3a.3.2.1) and $\g q=\g q_M$ from (1.3.2). The Hodge-Pink weights (see [HJ], below Remark 2.4) for this $\g q$: $\omega_1, \omega_2,\dots, \omega_r$ are 
$$-\g m, \dots, -\g m,\ \  -\g m+1, \dots, -\g m+1,\ \  \dots \dots,\ \ 0,\dots,0$$
where the number $-i$ occurs $k_{i+1}$ times. The Hodge-Pink filtration spaces $F^i\subset L(M)\otimes_{\n F_q[T]}\p[[N]]$ are the following: 
\medskip
$F^{-\g m}$ is the whole $L(M)\otimes_{\n F_q[T]}\p[[N]]$; 
\medskip
$\dim F^{-i}=k_1+k_2+...+k_{i+1}$ ($i=0,\dots, \g m$); $F^1=0$;
\medskip
$\dim Gr^{-i}_F=k_{i+1}$. 
\medskip
{\bf 3.18.} Let us continue here the exposition interrupted on the end of Section 3a.5. 
\medskip
Applying powers of $N$ to (3a.5.4.1), for any $v\in [0,\dots, \g m-1]$, $u\in [1,\dots, \g m+1]$ we can represent $N^{v}(\hat l_{u})$ as a linear combination of $N^z(\hat l_{y})$ where for a fixed $v$ the numbers $z, \ y$ satisfy $$z\in [v,\dots,\g m-1], \ \ y\in [z+2,\dots,\g m+1], \ \ \eqno{(3.19)}$$ Namely, there exist polynomials in $S_{****}$ denoted by $P_{uvyz}$ such that (matrix notations)

$$N^{v}\hat l_{u}=-\sum_{z=v}^{\g m-1}\sum_{y=z+2}^{\g m+1}P_{uvyz}N^z\hat l_{y}\eqno{(3.20)}$$
Clearly for $v=u-1$ we have $P_{uvyz}=S_{uvyz}$.
\medskip
{\bf 3.21.} The domain $v\ge u-1 \ \ \wedge \ \ \{z, \ y$ satisfy (3.19)\} is called the non-trivial domain of definition of $P_{****}$.
\medskip
For $v< u-1$ (trivial domain) we have:

$$ P_{u,v,y,z}=-I_*, \hbox{ resp. } P_{u,v,y,z}=0  \eqno{(3.22)}$$ for $y$, $z$ satisfying (3.19), $(y,z)=(u,v)$, resp. $(y,z)\ne(u,v)$.

\medskip

{\bf 3.23.} Example for $\g m=3$:

\medskip

$N^2\hat l_2=(S_{2131}S_{3242}-S_{2141})N^2\hat l_4$, i.e. $P_{2242}=-S_{2131}S_{3242}+S_{2141}$;

\medskip

$N^2\hat l_1=(-S_{1020}S_{2131}S_{3242}+S_{1020}S_{2141} +S_{1030}S_{3242}-S_{1040})N^2\hat l_4$, i.e.

\medskip

$P_{1242}=S_{1020}S_{2131}S_{3242}-S_{1020}S_{2141} -S_{1030}S_{3242}+S_{1040}$;

\medskip

$N\hat l_1=(S_{1020}S_{2131}-S_{1030})N\hat l_3+(S_{1020}S_{2141}-S_{1040})N\hat l_4+$

\medskip

$+(S_{1020}S_{2142}+S_{1031}S_{3242}-S_{1041})N^2\hat l_4$, i.e.

\medskip

$P_{1131}=-S_{1020}S_{2131}+S_{1030}$, $P_{1141}=-S_{1020}S_{2141}+S_{1040}$,

\medskip

$P_{1142}=-S_{1020}S_{2142}-S_{1031}S_{3242}+S_{1041}$.
\medskip
{\bf Remark 3.24. A.} Matrices $P_{****}$ are used to form a Siegel object of $M'$. In fact, not all $P_{uvyz}$ form it, but only $P_{u,v,y,y-2}$, see Definition 3.27 below. Hence, the set of these "essential" $P_{****}$ is also a tetrahedron.
\medskip
{\bf B.} Although we do not need this fact, let us give a formula for some $P_{****}$. Let us define a block unitriangular matrix $\g S$ whose $(i,j)$-th block is $S_{i,i-1,j,i-1}$ for $j>i$, $I_{k_i}$ for $i=j$ and 0 for $j<i$. Further on, we define a block unitriangular matrix $\g P$ whose $(i,j)$-th block is $-P_{i,j-2,j,j-2}$ for $j>i$, $I_{k_i}$ for $i=j$ and 0 for $j<i$. We have $\g P=\g S^{-1}$ (a proof follows immediately from the lemmas below).
\medskip
Some $P_{****}$ that enter in the below formula for $\bar B$ are not of the form of the elements of the inverse unitriangular matrix, for example $P_{1142}$, $\g m=3$.
\medskip
For the proof of Lemmas 3.32, 3.39, we need
\medskip
{\bf Lemma 3.25.} For all $i,\ j,\ \psi,\ \xi$ satisfying $i\in [2,\dots,\g m+1], \ j\in [1,\dots,\g m-i+2], \ \xi\in [\g m-j,\dots,\g m-1], \ \psi\in[\xi+2,\dots, \g m+1]$ we have
$$(\sum_{\be=0}^{j+\xi-\g m} \sum_{\al=i+\be}^{\g m+1-j+\be} S_{i-1,i-2,\al,i-2+\be} P_{\al,\g m-j+\be,\psi,\xi})-$$ $$-S_{i-1,i-2,\psi,i-2+\xi+j-\g m}+P_{i-1,\g m-j,\psi,\xi}=0\eqno{(3.25.1)}$$
(A recurrent formula for $P_{****}$).
\medskip
{\bf Proof.} First, we rewrite (3.16) for $u=i-1$:

$$N^{i-2}\hat l_{i-1}=-\sum_{z=i-2}^{\g m-1}\sum_{y=z+2}^{\g m+1}S_{i-1,i-2,y,z}N^z\hat l_{y}\eqno{(3.25.2)}$$

Now, for any $$j=1,\dots,\g m-i+2\eqno{(3.25.3)}$$ we apply $N^{\g m-i+2-j}$ to (3.25.2):

$$N^{\g m-j}\hat l_{i-1}=-\sum_{z=i-2}^{i+j-3}\sum_{y=z+2}^{\g m+1}S_{i-1,i-2,y,z}N^{z+\g m-i+2-j}\hat l_{y}\eqno{(3.25.4)}$$

(since $N^\g m=0$, we get that $z\le i+j-3$ ).

We change the summation variables: $z \to i-2+\be$, and $y \to \al$, we get

$$N^{\g m-j}\hat l_{i-1}=-\sum_{\be=0}^{j-1}\sum_{\al=i+\be}^{\g m+1}S_{i-1,i-2,\al,i-2+\be}N^{\g m+\be-j}\hat l_{\al}\eqno{(3.25.5)}$$

Now we use (3.20) making the following variable change:

$$u \to \al\ \ \ \ \ y\to \psi \ \ \ \ \ v\to \g m-j+\be\ \ \ \ \ z\to \xi.$$
We get
$$N^{\g m-j+\be}\hat l_\al=-\sum_{\xi=\g m-j+\be}^{\g m-1} \sum_{\psi=\xi+2}^{\g m+1} P_{\al,\g m-j+\be,\psi,\xi}N^{\xi}\hat l_\psi\eqno{(3.20a)}$$

We substitute (3.20a) in (3.25.5):

$$N^{\g m-j}\hat l_{i-1}=\sum_{\be=0}^{j-1}\sum_{\al=i+\be}^{\g m+1} \sum_{\xi=\g m-j+\be}^{\g m-1}\sum_{\psi=\xi+2}^{\g m+1} S_{i-1,i-2,\al,i-2+\be} P_{\al,\g m-j+\be,\psi,\xi}N^{\xi}\hat l_\psi\eqno{(3.25.6)}$$

We change the order of summation in (3.25.6):

$$N^{\g m-j}\hat l_{i-1}=\sum_{\xi=\g m-j}^{\g m-1} \sum_{\psi=\xi+2}^{\g m+1} (\sum_{\be=0}^{j+\xi-\g m} \sum_{\al=i+\be}^{\g m+1} S_{i-1,i-2,\al,i-2+\be} P_{\al,\g m-j+\be,\psi,\xi}) N^{\xi}\hat l_\psi\eqno{(3.25.7)}$$

We rewrite (3.20) making changes:

$$u\to i-1\ \ \ \ \ \ y \to \psi\ \ \ \ \ \ v\to \g m-j\ \ \ \ \ \  z\to \xi.$$
We get
$$N^{\g m-j}\hat l_{i-1}=-\sum_{\xi=\g m-j}^{\g m-1}\sum_{\psi=\xi+2}^{\g m+1} P_{i-1,\g m-j,\psi,\xi}N^\xi \hat l_{\psi}\eqno{(3.25.8)}$$

For $\psi\ge \xi+2$ elements $N^\xi l_{\psi i}$, $i=1,\dots,k_\psi$, are linearly independent over $\p$. Hence, (3.25.7), (3.25.8) imply

$$P_{i-1,\g m-j,\psi,\xi}=-\sum_{\be=0}^{j+\xi-\g m} \sum_{\al=i+\be}^{\g m+1} S_{i-1,i-2,\al,i-2+\be} P_{\al,\g m-j+\be,\psi,\xi}\eqno{(3.25.9)}$$

Here the domain of $\xi$, $\psi$ is:

$$\xi\in [\g m-j,\dots,\g m-1], \ \ \ \psi\in[\xi+2,\dots, \g m+1]$$
Taking into consideration (3.22) we can rewrite (3.25.9) as follows:

$$P_{i-1,\g m-j,\psi,\xi}=-(\sum_{\be=0}^{j+\xi-\g m} \sum_{\al=i+\be}^{\g m+1-j+\be} S_{i-1,i-2,\al,i-2+\be} P_{\al,\g m-j+\be,\psi,\xi})+$$ $$+S_{i-1,i-2,\psi,i-2+\xi+j-\g m}\eqno{(3.25.10)}$$ with the same domain of $\xi$, $\psi$. This is (3.25.1). Because of (3.25.3), this formula is valid for $\g m-j\ge i-2$ (the non-trivial case of the definition of $P_{****}$). $\square$
\medskip
Let us consider the symmetry $\g s: \n Z^4\to\n Z^4$ defined as follows: $\g s(\al,\be,\ga,\de)=(\g m+2-\ga,\ \g m-1-\de,\ \g m+2-\al,\ \g m-1-\be)$.
\medskip
{\bf Remark 3.26.} $\g s$ has the following geometric interpretation. Let us consider a matrix $NL$ whose $(i,j)$-th entry is a symbol $N^{i-1}\hat l_{j}$. We interpret a quadruple $(\al,\be,\ga,\de)$ as a vector from $N^\be \hat l_\al$ to $N^\de \hat l_\ga$ in $NL$. $\g s$ is the reflection of this vector with respect to the center of $NL$ and the inversion of its direction.
\medskip
{\bf Definition 3.27.} $\bar S_{uvyz}:=-P_{\g s(uvyz)}^t$ (defined if $P_{\g s(uvyz)}$ has meaning).
\medskip
We consider block $r\times r$-matrices having the following block structure: their block size is $(\g m+1)\times(\g m+1)$, quantities of columns in blocks are $k_{\g m+1},k_\g m,\dots,k_1$ (counting from the left to the right), and quantities of lines in blocks are $k_{1},k_2,\dots,k_{\g m+1}$ (counting from up to down). Hence, the $(\al,\be)$-th block of this matrix is a $k_{\al}\times k_{\g m+2-\be}$-matrix. These matrices will be called skew $k_*$-block matrices.

%We denote their sequence $(k'_1,\dots,k'_{\g m+1})$ of dual numbers by $k'_*$.
\medskip
$\forall \ i = 0,\dots, \g m$ we define skew $k_*$-block matrices $C_i=C_i(S_{****})$ as follows:
\medskip
The $(\al, \be)$-th block of $C_i$ is $S^t_{\g m+2-\be,\g m+1-\be,\al,i}$ if the quadruple $(\g m+2-\be,\g m+1-\be,\al,i)$ satisfies (3a.5.2, 3a.5.3)\footnotemark \footnotetext{Obviously it always satisfies (3a.5.3).} (i.e. if $S_{\g m+2-\be,\g m+1-\be,\al,i}$ exists);
the $(i+1,\g m+1-i)$-th block of $C_i$ is $I_{k_{i+1}}$, all other blocks of $C_i$ are 0. Namely,
$$(C_i)_{\al\be}=S^t_{\g m+2-\be,\g m+1-\be,\al,i}\eqno{(3.28.1)}$$
$$(C_i)_{i+1,\g m+1-i}=I_{k_{i+1}}\eqno{(3.28.2)}$$

$\forall \ i = 0,\dots, \g m$ we define skew $k_*$-block matrices $\bar C_i=\bar C_i(S_{****})$ as follows:
\medskip
The $(\al,\be)$-th block of $\bar C_i$ is given by the formula

$$(\bar C_i)_{\al,\be}=-P_{\al,\g m-1-i,\g m+2-\be,\g m-\be}=\bar S^{t}_{\be,\be-1,\g m+2-\al,i}\eqno{(3.29)}$$
if the quadruple $(\al,\g m-1-i,\g m+2-\be,\g m-\be)$ belongs to the non-trivial domain of $P_{****}$;
\medskip
For $i=0,\dots,\g m$
$$(\bar C_i)_{\g m+1-i,i+1}=I_{k_{\g m+1-i}}\eqno{(3.30)}$$
other block entries of $\bar C_i$ are 0.
\medskip
{\bf Remark 3.31.} Formula (3.30) is concordant with (3.29), if we consider $P_{****}$ from (3.22). Nevertheless, some 0-blocks of $\bar C_i$ correspond to $P_{**yz}$ where $(y,z)$ do not satisfy (3.19), and hence this $P_{****}$ is not defined.
\medskip
Finally, we define elements $B(S_{****}):=\sum_{i=0}^\g m C_iN^i\in M_r(\p)[N]$ and $\bar B(S_{****}):=\sum_{i=0}^\g m \bar C_iN^i\in M_r(\p)[N]$.
\medskip
Example for $\g m=3$:
$$B(S_{****})=\left(\matrix 0&0&0&I_{k_1}\\0&0&0& S^t_{1020}\\0&0&0& S^t_{1030}\\0&0&0& S^t_{1040} \endmatrix \right)+
\left(\matrix 0&0&0&0\\0&0&I_{k_2}& 0\\0&0&S^t_{2131}& S^t_{1031}\\0&0&S^t_{2141}& S^t_{1041} \endmatrix \right)N+$$ $$+ \left(\matrix 0&0&0&0\\0&0&0& 0\\0&I_{k_3}&0& 0\\0& S^t_{3242}&S^t_{2142}& S^t_{1042} \endmatrix \right)N^2 + \left(\matrix 0&0&0&0\\0&0&0&0\\0&0&0& 0\\I_{k_4}&0&0&0 \endmatrix \right)N^3$$
\medskip
$$\bar B(S_{****})=\left(\matrix \bar S^{t}_{1040}&0&0&0&\\ \bar S^{t}_{1030}&0&0&0\\ \bar S^{t}_{1020}&0&0&0\\I_{k_4}&0&0&0  \endmatrix \right)+
\left(\matrix \bar S^{t}_{1041}&\bar S^{t}_{2141}&0&0\\ \bar S^{t}_{1031}&\bar S^{t}_{2131}&0& 0\\0&I_{k_3}&0&0 \\0&0&0&0  \endmatrix \right)N+$$ $$+ \left(\matrix \bar S^{t}_{1042}&\bar S^{t}_{2142}& \bar S^{t}_{3242}&0\\0&0&I_{k_2}& 0\\0&0&0& 0\\0&0&0& 0  \endmatrix \right)N^2 + \left(\matrix 0&0&0&I_{k_1}\\0&0&0&0\\0&0&0& 0\\0&0&0&0 \endmatrix \right)N^3=$$
\medskip
$$=\left(\matrix -P_{1242}&0&0&0&\\-P_{2242}&0&0&0\\-P_{3242}&0&0&0\\I_{k_4}&0&0&0  \endmatrix \right)+
\left(\matrix -P_{1142}&-P_{1131}&0&0\\-P_{2142}&-P_{2131}&0& 0\\0&I_{k_3}&0&0 \\0&0&0&0  \endmatrix \right)N+$$ $$+ \left(\matrix -P_{1042}&-P_{1031}& -P_{1020}&0\\0&0&I_{k_2}& 0\\0&0&0& 0\\0&0&0& 0  \endmatrix \right)N^2 + \left(\matrix 0&0&0&I_{k_1}\\0&0&0&0\\0&0&0& 0\\0&0&0&0 \endmatrix \right)N^3$$
\medskip
{\bf Lemma 3.32.} $B(S_{****})^t \cdot \bar B(S_{****})=I_r N^\g m\in M_r(\p)[N]$.
\medskip
{\bf Proof.} We denote $B(S_{****})^t \cdot \bar B(S_{****})$ by $\sum_\mu \Cal C_\mu N^\mu$. The fact that $\Cal C_\g m=I_{r}$ is obvious: the only non-zero factors that enter in the sum $\sum_{\ga=0}^{\g m} C_\ga^t \bar C_{\g m-\ga}$ are products of blocks of $C_*$, $\bar C_*$ containing $I_*$, and they form $I_r$. Also it is obvious that for $\mu>\g m$ we have $\Cal C_\mu=0$, because all products whose sum is $\Cal C_\mu$, have at least one factor 0. We need to consider $\Cal C_\mu$ for $\mu<\g m$. (3.28.1), (3.28.2), (3.29), (3.30) give us (here and below $(C^t_\ga)_{\nu\de}$ is the $(\nu\de)$-th block of $C^t_\ga$, i.e. $(C^t_\ga)_{\nu\de}=((C_\ga)_{\de\nu})^t$ )
$$(\Cal C_\mu)_{\nu\pi}=\sum_{\ga=0}^\mu \sum_{\de=1}^{\g m+1} (C^t_\ga)_{\nu\de}(\bar C_{\mu-\ga})_{\de\pi}\eqno{(3.32.1.1)}$$ $$=-\sum_{\ga,\de} S_{\g m+2-\nu, \g m+1-\nu, \de, \ga}P_{\de,\g m-1-\mu+\ga,\g m+2-\pi,\g m-\pi}+ \eqno{(3.32.1.2)}$$
$$-P_{\g m+2-\nu,2\g m-\mu-\nu,\g m+2-\pi,\g m-\pi}+S_{\g m+2-\nu,\g m+1-\nu,\g m+2-\pi,\mu+1-\pi}\eqno{(3.32.1.3)}$$
where (3.32.1.2) corresponds to the products $(C^t_\ga)_{\nu\de}(\bar C_{\mu-\ga})_{\de\pi}$ where both terms $\ne 0, \ I_*$, and (3.32.1.3) corresponds to the products where one of the terms is $I_*$.

Let us find the relations satisfied by $\mu,\ \nu, \ \pi$ and the domain of summation by $\ga, \ \de$ in (3.32.1.2). We have
$$\matrix (C^t_\ga)_{\nu\de}\ne0, I_{k_*} \ \iff \ \nu\ge \g m+1-\ga \ \ \wedge \ \ \de\ge\ga+2 \\
(\bar C_{\mu-\ga})_{\de\pi}\ne0, I_{k_*} \ \iff \ \de\le \g m-(\mu-\ga) \ \ \wedge \ \ \pi \le \mu-\ga+1\endmatrix \eqno{(3.32.2)}$$
The set of $\ga, \ \de$ is non-empty $\iff \ \mu\le \g m-2$ and $\mu +\nu -\pi\ge \g m$. In this case the conditions (3.32.2) on $\ga, \ \de$ become
$$\mu+1-\pi\ge\ga\ge \g m+1-\nu\eqno{(3.32.3.1)}$$
$$\g m+\ga-\mu\ge\de\ge\ga+2\eqno{(3.32.3.2)}$$
Now we use Proposition 3.25 for
$$\matrix i=\g m+3-\nu && \mu=j+i-3 \\
j=\mu+\nu-\g m &\iff & \nu=\g m-i+3 \\
\psi=\g m-\pi+2 && \pi=\g m-\xi=\g m-\psi+2\\
\xi=\g m-\pi \endmatrix\eqno{(3.32.4)}$$
and summation variables $\al, \ \be$ in (3.25.1) are
$$\matrix \al=\de \\ \be=\ga-i+2 \endmatrix\eqno{(3.32.5)}$$
Under this variable change, (3.32.1.2) becomes the double sum in (3.25.10), and (3.32.1.3) becomes
$$-P_{i-1,\g m-j,\psi,\xi}+S_{i-1,i-2,\psi,i-2+\xi+j-\g m}$$
hence the desired. $\square$
\medskip
Let for $i=0,\dots,\g m-1$ $X_i$ be skew $k_*$-matrices having the following property:
\medskip
If $(\al, \ \be)$ are such that the $(\al,\ \be)$-block of $\bar C_i$ is 0 or $I_*$ then $(X_i)_{\al\be}=(\bar C_i)_{\al\be}$;

If $(\al, \ \be)$ are such that the $(\al,\ \be)$-block of $\bar C_i$ is $\ne 0, \ I_*$ then $(X_i)_{\al\be}$ is arbitrary.
\medskip
We denote $X:=\sum_{i=0}^{\g m-1}X_iN^i$.
\medskip
{\bf Lemma 3.33.} If $B(S_{****})^t \cdot X\in N^\g mM_r(\p[N])$ then $X=\bar B(S_{****})$.
\medskip
{\bf Proof.} For any fixed $\mu, \ \nu, \ \pi$ (3.32.1.2), (3.32.1.3) become
$$\sum_{\ga,\de} S_{\g m+2-\nu, \g m+1-\nu, \de, \ga} (X_{\mu-\ga})_{\de,\pi}\eqno{(3.33.1a)}$$ $$
+ (X_{\mu+\nu-\g m-1})_{\g m+2-\nu,\pi}+S_{\g m+2-\nu,\g m+1-\nu,\g m+2-\pi,\mu+1-\pi}=0\eqno{(3.33.1b)}$$
This is system of linear equations with unknowns $(X_i)_{\al\be}$ where $$ 0\le i \le \g m-1 \ \ \ \ 1\le \be\le i+1\ \ \ \ 1\le\al \le \g m-i \eqno{(3.33.2)}$$ (for other values of $i,\ \al, \ \be$ we have $(X_i)_{\al\be}=0$ or $I_*$). We arrange $(X_i)_{\al\be}$ in decreasing order of $i+\al$ (for $(X_i)_{\al\be}$ having equal $i+\al$ their ordering is arbitrary), and we arrange equations (3.33.1) in decreasing order of $\mu$ (for equations having equal $\mu$ the order of equations corresponds to the order of $(X_i)_{\al\be}$ having equal $i+\al$). Under this arrangement of unknowns and equations, the matrix of the system (3.33.1) becomes unitriangular. Really, for any $i,\ \al, \ \be$ satisfying (3.33.2) there is exactly one values of $\mu, \ \nu, \ \pi$ --- namely, $$ \mu=i+\al-1 \ \ \ \ \nu=\g m+2-\al\ \ \ \ \pi=\be  $$ such that the first term of (3.33.1b) is $(X_i)_{\al\be}$. Other terms of (3.33.1) for these $\mu, \ \nu, \ \pi$ --- namely, the terms that enter in (3.33.1a) --- contain $(X_{\mu-\ga})_{\de,\pi}$ such that $\mu-\ga+\de>i+\al$ hence unitriangularity.
\medskip
Lemma 3.32 affirms that
$$(X_i)_{\de,\pi}=-P_{\de,\g m-1-i,\g m+2-\pi,\g m-\pi}$$
is a solution to this system. Unitriangularity implies that this solution is unique. $\square$
\medskip
Let us consider (1.3.2) for the present setting. We have (notations of (0.1)) $R=\n F_q[\th\cdot I_n+N]$. (1.3.2) becomes
$$0\to\g q_L\to L\underset{R}\to{\otimes}\p[[N]]\to V\to0\eqno{(3.33a)}$$
By abuse of notations, we shall denote elements $l_i\in L\subset V$ as elements of $V$, and elements $l_i=l_i\otimes 1\in L\underset{R}\to{\otimes}\p[[N]]$, by the same symbol: there will be no confusion. For example, $N^{\g m}l_i\ne0$ in $L\underset{R}\to{\otimes}\p[[N]]$ while $N^{\g m}l_i=0$ in $V$.
\medskip
{\bf Lemma 3.34.} $\forall \ u=1,\dots, \g m+1$, for $v=u-1$, $\forall \ i=1,\dots, k_u$ the elements

$$\om_{ui}:=N^{v}l_{ui}+\sum_{y=u+1}^{\g m+1}\sum_{z=u-1}^{y-2}\sum_{j=1}^{k_y}(S_{uvyz})_{ij}N^zl_{yj}\eqno{(3.34.1)}$$

form a basis of $\g q_L$. $\square$
\medskip
We denote the set of elements $\om_{ui}$ ($u$ is fixed, $i$ varies) by $\hat \om_{u}$ (matrix columns). So, (3.34.1) becomes ($v=u-1$)
$$\hat \om_{u}=N^{v}\hat l_{u}+\sum_{z=u-1}^{\g m-1}\sum_{y=z+2}^{\g m+1}S_{uvyz}N^z\hat l_{y}\eqno{(3.35)}$$
\medskip
{\bf 3.35.1.} Let $L'=\Hom_R(L,R)$ be the dual module, and let $\la_i$, $i=1,...,r$ be the basis of $L'$ dual to $\hat l$, i.e. $\la_i(l_j)=\delta_{ij}$. We shall need the dual numbers $k'_i:=k_{\g m+2-i}$ (inverse order of $k_*$). We consider the analogous two-subscript notation of $\la_i$, but the order of segments of the partition of $\la_i$ is opposite, namely:
$$(\la_{\g m+1,1},\dots,\la_{\g m+1,k'_{\g m+1}},\ \ \la_{\g m,1},\dots,\la_{\g m,k'_{\g m}},\ \ \ \dots \ \ \ , \la_{11},\dots,\la_{1,k'_{1}}):= (\la_1,\dots, \la_r)\eqno{(3.35.2)}$$
(order of elements $\la_*$ is the same in both sides of this equality). Further, let $\g q'_L$ be defined by the same formula (1.6), $L'$ instead of $L(M)'$.
\medskip
{\bf Lemma 3.36.} $\forall \ u=1,\dots, \g m+1$, for $v=u-1$, $ \forall \ i=1,\dots, k'_u$ the elements

$$\chi_{ui}:=N^{v}\la_{ui}+\sum_{y=u+1}^{\g m+1}\sum_{z=u-1}^{y-2}\sum_{j=1}^{k'_y}(\bar S_{uvyz})_{ij}N^z\la_{yj}\eqno{(3.36.1)}$$

form a basis of $\g q'_L$.
\medskip
{\bf Proof.} As above we denote the set of elements $\la_{ui}$, resp. $\chi_{ui}$ ($u$ is fixed, $i$ varies) by $\hat \la_{u}$, resp. $\hat \chi_{u}$ (matrix columns). (3.35), (3.36.1) can be written in terms of blocks of $C_i$, $\bar C_i$:

$$\hat \om_{u}=\sum_{z=0}^{\g m}\sum_{y=1}^{\g m+1} (C_z^t)_{\g m+2-u,y}N^z\hat l_{y}$$
$$\hat \chi_{u}=\sum_{z=0}^{\g m}\sum_{y=1}^{\g m+1} (\bar C_z^t)_{uy}N^z\hat \la_{\g m+2-y}$$
We must prove that $\forall \ u_1, \ u_2$ we have $\hat \om_{u_1}\hat \chi_{u_2}^t=\de_{u_1}^{u_2}I_{k_{u_1}}N^\g m$ (product is pairing). This is immediate:

$$\hat \om_{u_1}\hat \chi_{u_2}^t=\sum_{z_1=0}^{\g m}\sum_{y_1=1}^{\g m+1}\sum_{z_2=0}^{\g m}\sum_{y_2=1}^{\g m+1} (C_{z_1}^t)_{\g m+2-u_1,y_1}N^{z_1}\hat l_{y_1} N^{z_2}\hat \la_{\g m+2-y_2}^t (\bar C_{z_2})_{y_2,u_2}$$
We have $\hat l_{y_1} \hat \la_{\g m+2-y_2}^t=\de_{y_1}^{y_2}I_{k_{y_1}}$, hence
$$\hat \om_{u_1}\hat \chi_{u_2}^t=\sum_{z_1=0}^{\g m}\sum_{y=1}^{\g m+1}\sum_{z_2=0}^{\g m} (C_{z_1}^t)_{\g m+2-u_1,y}N^{z_1+z_2} (\bar C_{z_2})_{y,u_2} = \sum_{z=0}^{2\g m}(\Cal C_z)_{\g m+2-u_1,u_2}N^z$$
Lemma 3.32 implies the desired. $\square$
\medskip
{\bf Corollary 3.37.} Matrices $S_{uvyz}(L')$ for the dual lattice $L'$ are $\bar S_{uvyz}(L)$ (order of segments of $\la_*$, and hence of numbers $k_*$, is inverse).
\medskip
{\bf 3.37A.} Until now we considered abstract $L$ and $V$. Now we consider $L$, $V$ coming from $M$. Namely, let $M$ be an uniformizable t-motive such that $N^\g m=0$. We let $L=L(M)$, $V=\Lie(M)$, and we consider $\hat l_i$, $S_*$ etc. for them. (3a.1.3.1) holds for this case, because $L(M)\underset{\n F_q[T]}\to{\otimes}\p[[T-\th]]\to \Lie(M)$ is an epimorphism, see (1.3.2).

Let $$\Psi_N=\sum_{u=-\g m}^{\infty}D_{-u}N^u $$ be the $\th$-shift of $\Psi(M)$ (the fact that the series $\Psi_N$ exists and starts from $D_\g mN^{-\g m}$ follows from Proposition 2.4). We represent each $D_i$ ($i=1,\dots, \g m$) as a union of $r\times k_j$-blocks $D_{ij}$, $j=1,\dots, \g m+1$, namely $D_{i1}$ is a submatrix of $D_i$ formed by its first $k_1$ columns, $D_{i2}$ is a submatrix of $D_i$ formed by its next $k_2$ columns, etc, until $D_{i,\g m+1}$ is a submatrix of $D_i$ formed by its last $k_{\g m+1}$ columns.
\medskip
Proposition 2.4 implies that $(D_{ij})_{\al\be}=-<N^{i-1}(l_{j\be}),f_\al>$. This means that there are relations between $D_{ij}$ coming from (3.20), namely:

$$D_{v+1,u}= -\sum_{z=v}^{\g m-1}\sum_{y=z+2}^{\g m+1}D_{z+1,y} P_{uvyz}^t\eqno{(3.38)}$$

Like in (3.21), for $(z,y)$ satisfying $y\ge z+1$ (resp. $y< z+1$) we shall call the corresponding $D_{zy}$ as belonging to the trivial (resp. non-trivial) domain.
\medskip
{\bf Lemma 3.39.} $\Psi_N B(S_{****})\in M_r(\p[[N]])$.
\medskip
{\bf Proof.} For any $\mu=0,\dots, \g m-1$ we must prove that $\g E_\mu:=\sum_{\de =0}^\mu D_{\g m-\de}C_{\mu-\de}$ is 0. The $\nu$-th block ($\nu=1,\dots,\g m+1$) of this matrix is
$$(\g E_\mu)_\nu:=\sum_{\de =0}^\mu \sum_{\ga=1}^{\g m+1} D_{\g m-\de,\ga}(C_{\mu-\de})_{\ga \nu}\eqno{(3.39.1)}$$
We have
$$(C_{\mu-\de})_{\ga \nu}\ne0, I_{k_*} \ \ \iff \de\le \nu+\mu-\g m-1 \ \ \wedge \ \ \ga\ge \mu-\de+2$$
$$(C_{\mu-\de})_{\ga \nu}=I_{k_*} \ \ \iff \de=\nu+\mu-\g m-1 \ \ \wedge \ \ \ga=\g m+2-\nu$$
hence (3.39.1) becomes
$$(\g E_\mu)_\nu=\sum_{\de =0}^{\nu+\mu-\g m-1}\sum_{\ga=\mu-\de+2}^{\g m+1}D_{\g m-\de,\ga}S^t_{\g m+2-\nu,\g m+1-\nu,\ga,\mu-\de}+\eqno{(3.39.2.1)}$$
$$+D_{2\g m+1-\nu-\mu,\g m+2-\nu}\eqno{(3.39.2.2)}$$
where (3.39.2.1) is non-empty if $\nu+\mu\ge \g m+1$.

Terms $D_{2\g m+1-\nu-\mu,\g m+2-\nu}$ always belong to the non-trivial domain. We separate the terms of (3.39.2.1) in terms of trivial and non-trivial domain:
$$(\g E_\mu)_\nu=\sum_{\de =0}^{\nu+\mu-\g m-1}\sum_{\ga=\mu-\de+2}^{\g m-\de} D_{\g m-\de,\ga}  S^t_{\g m+2-\nu,\g m+1-\nu,\ga,\mu-\de}+\eqno{(3.39.3.1)}$$
$$+\sum_{\de =0}^{\nu+\mu-\g m-1}\sum_{\ga=\g m-\de+1}^{\g m+1}  D_{\g m-\de,\ga}  S^t_{\g m+2-\nu,\g m+1-\nu,\ga,\mu-\de}+\eqno{(3.39.3.2)}$$
$$+D_{2\g m+1-\nu-\mu,\g m+2-\nu}\eqno{(3.39.3.3)}$$

Now we substitute non-trivial $D_{**}$ by linear combinations of the trivial ones, using (3.38):
$$(\g E_\mu)_\nu=-\sum_{\de =0}^{\nu+\mu-\g m-1}\sum_{\ga=\mu-\de+2}^{\g m-\de}   \sum_{z=\g m-\de-1}^{\g m-1}\sum_{y=z+2}^{\g m+1}D_{z+1,y} P_{\ga,\g m-\de-1,y,z}^t       S^t_{\g m+2-\nu,\g m+1-\nu,\ga,\mu-\de}+\eqno{(3.39.4.1)}$$
$$+\sum_{\de =0}^{\nu+\mu-\g m-1}\sum_{\ga=\g m-\de+1}^{\g m+1}D_{\g m-\de,\ga}S^t_{\g m+2-\nu,\g m+1-\nu,\ga,\mu-\de}-\eqno{(3.39.4.2)}$$
$$-\sum_{z=2\g m-\nu-\mu}^{\g m-1}\sum_{y=z+2}^{\g m+1}D_{z+1,y} P_{\g m+2-\nu,2\g m-\nu-\mu,y,z}^t \eqno{(3.39.4.3)}$$
Now we change variables in (3.39.4.2):
$$\de=\g m-z-1$$
$$\ga=y$$
interchange the order of summation and transpose:
$$(\g E_\mu)_\nu=\sum_{z=2\g m-\nu-\mu}^{\g m-1}\sum_{y=z+2}^{\g m+1} \g K(\mu,\nu,z,y) D^t_{z+1,y}$$
where $$\g K(\mu,\nu,z,y)=-( \sum_{\de =\g m-z-1}^{\nu+\mu-\g m-1}\sum_{\ga=\mu-\de+2}^{\g m-\de} S_{\g m+2-\nu,\g m+1-\nu,\ga,\mu-\de}    P_{\ga,\g m-\de-1,y,z}) +$$
$$+S_{\g m+2-\nu,\g m+1-\nu,y,\mu-\g m+z+1}-P_{\g m+2-\nu,2\g m-\nu-\mu,y,z} \eqno{(3.39.5)}$$
Change of variables in (3.39.5):
$$\matrix y=\psi & \nu=\g m+3-i & \ga=\al \\ z=\xi & \mu=i+\al-3 & \de=\al-\be-1 \endmatrix $$
transforms (3.39.5) to (3.25.1), hence all $\g K(\mu,\nu,z,y)$ are 0. $\square$
\medskip
\medskip
{\bf 4. Proof.}
\medskip
Let $M$ be such that $M'=M^{\prime \g m}$ --- its $\g m$-dual --- exists and satisfies $N^{\g m}=0$, where $N$ is of $M'$. We denote numbers $k_i$, resp. $k'_i$, for $L(M)$, by $k_i(M)$, resp. $k'_i(M)$.
\medskip
{\bf Lemma 4.1.} For $i=1,\dots, \g m+1$ we have $k_i(M')=k'_{i}(M)$.
\medskip
{\bf Proof.} For $\g m=1$ this is [GL07], Lemma 10.2. For any $\g m$ the proof is analogous, let us give it. We denote numbers $d_i$, $c_i$ (see (3.4.1) and below) for $M'$ by $d'_i$, $c'_i$ respectively. We need an explicit formula for the action of $\tau$ on $M'$. Namely, there exists a matrix $Q=Q(M,\hat f)\in M_r(\p[T])$ such that $Q\hat f=\tau \hat f$. The matrix $Q$ defines $M$ uniquely. We denote $Q(M', \hat f')$ by $Q'$. We have (this follows immediately from the definition of $M'$; see also [GL07], (1.10.1)):
$$Q'=(T-\th)^\g mQ^{t-1}$$
(here and below $*^{t-1}$ means $(*^t)^{-1}$).

The module $\tau M$ is a $\p[T]$-submodule of $M$ (because $a \tau x =
\tau a^{1/q} x$ for
$x\in M$),
hence there are $\p[T]$-bases $f_*=(f_1, ..., f_r)^t$, $g_*=(g_1, ...,
g_r)^t$ of $M$,
$\tau M$ respectively such that $g_i=P_if_i$, where $P_r | P_{r-1} | ... |
P_1$, \ \ $P_i \in
\p[T]$. Condition (1.2.1) means that $\forall i$ \ $(T-\theta)^\g m f_i
\in \tau M$, i.e.
$P_i|(T-\theta)^\g m$. There exists an isomorphism of $\p[T]$-modules $\Lie(M)$ and $M/\tau M$, hence $\forall i$ \ $P_i=(T-\theta)^{d_i}$, where $d_i$ are from (3a.3.2) (numbers $d_i$ are $m_{r+1-i}$ from [GL07], Lemma 10.2). There exists a matrix $\goth Q=\{\goth q_{ij}\}\in
M_r(\p[T])$ such that
$$\tau f_i = \sum_{j=1}^r \goth q_{ij}g_j= \sum_{j=1}^r \goth q_{ij} P_j
f_j\eqno{(4.1.1)}$$
Although $\tau$ is not a linear operator, it is easy to see that
$\goth Q\in GL_r(\p[T])$
(really, there exists $C=\{c_{ij}\}\in M_r(\p[T])$ such that $g_i=P_i
f_i=\tau (\sum_{j=1}^r
c_{ij}f_j)$, we have $C^{(1)}\goth Q =I_r$).

We denote the matrix $\diag(P_1, P_2, ... ,P_r)$ by $\Cal P$, so
(4.1.1) means that $Q= \goth Q \Cal P$. This implies that
$Q'=Q({M'})$ satisfies
$$Q'=\goth Q^{t-1}\diag((T-\theta)^{\g m-d_r}, ...,
(T-\theta)^{\g m-d_1})$$
This means that $\forall \ i=1,\dots,r$ \ $d'_i=\g m-d_{r+1-i}$ and hence $\forall \ i=1,\dots,\g m$ \ $c'_i=r-c_{\g m+1-i}$. The lemma follows from (3a.3.2.1). $\square$
\medskip
We have $L(M')=L(M)'$, see (2.1), (2.1a). Hence, we identify $\la_*$ from (3.35.1) with $\vf_*$ from (2.1a). We shall consider another basis $\hat \eta$ of $L(M')$ obtained by a permutation matrix from the basis (3.35.2). Namely, we let $\eta_{ij}:=\vf_{ij}$ ($i=1, \dots, \g m+1, \ j=1,\dots, k'_i$), but the order of elements $\eta_{ij}$ is the following ($\hat \eta$ is a matrix column):

$$\hat \eta=(\eta_{11},\dots, \eta_{1,k'_1}, \eta_{21},\dots, \eta_{2,k'_2}, \dots, \eta_{\g m+1,1},\dots, \eta_{\g m+1,k'_{\g m+1}})^t$$

We shall consider only $M$ satisfying the following (compare with 3a.4.6)
\medskip
{\bf Condition 4.2.} $M'$ exists, satisfies $N^{\g m}=0$, and for any $u$ the elements $N^\al(\eta_{\be\ga})$ (where $\al\in[u,\dots, \g m-1]$, $\be\in[\al+2,\dots, \g m+1]$, $\ga\in [1,\dots, k'_\be]$) are linearly independent over $\p$.
\medskip
According the general principle that almost all $n$-uples $(v_1,\dots,v_n)$ of vectors in $n$-dimensional vector space form a basis of this space, we can guess that almost all $M$ satisfy 4.2. Really, Lemma 4.1 affirms that the dimension of $N^u \Lie(M')$ is exactly the quantity of elements $N^\al(\eta_{\be\ga})$ mentioned in 4.2. Again by Lemma 4.1, we see that Condition 4.2 implies that these elements are a basis of $N^u \Lie(M')$.
\medskip
{\bf Remark 4.3.} We can use ideas of [GL17] in order to prove rigorously that a large class of Anderson t-motives $M$ satisfies Condition 4.2. Namely, the subject of [GL17] is consideration of Anderson t-motives $M(A)$ of dimension $n$ and rank $2n$ given by the equation $$Te_*=\theta e_*+A\tau e_* +\tau^2e_* \eqno{(4.3.1)}$$ ([GL17], (0.8)) where $A\in M_n(\p)$ is near 0. A Siegel matrix of $M(0)$ is $\om I_n$ where $\om\in \n F_{q^2}-\n F_q$. It is shown in [GL17] that any matrix near $\om I_n$ is a Siegel matrix of some $M(A)$.

We can apply methods of [GL17] to standard-1 t-motives (see [GL07], (11.3); (11.2.1)) having $N\ne0$. Let us give here its definition. Let $N=N_{ij}$ be the matrix of $N$, i.e. the Jordan matrix with 0-blocks of sizes $d_1,d_2,...,d_\al$ from (3a.1.1). We consider a function (see [GL07], (11.2): $\g k: (1, ..., n) \to \n Z^+$ where $\n Z^+$ is the set of integers $\ge 1$. A standard-1 t-motive $M(N, \g k, a_{***})$ is an Anderson
t-motive of dimension $n$ given by the formulas ($i=1,...,n$):
$$Te_{i}= (\theta e_{i}+\sum_{\alpha=1}^n N_{i\al}e_\al )
+\sum_{\alpha=1}^n\sum_{j=1}^{k(\alpha)-1}a_{j,i,\alpha}\
\tau^je_\alpha + \tau^{\g k(i)}e_i\eqno{(4.3.2)}$$ where $a_{j,i,\alpha}\in\p$ is the $(i,\alpha)$-th entry of the matrix $\goth A_j$ (see (1.3)). This is a simplified version of [GL07], (11.2.1). As an initial t-motive (analog of $M(0)$ of [GL17]) we choose $M(N, \g k, a_{***})$ for $a_{***}=0$. We denote the set of its $S_{****}$ by $S(N, \g k, 0)$.

Methods of [GL17] show that for any $S_{****}$ near $S(N, \g k, 0)$ there exist $a_{***}$ near 0 such that a Siegel set of $M(N, \g k, a_{***})$ is the given $S_{****}$. It is clear that for almost all $a_{***}$ near 0 the Condition 4.2 holds for $M(N, \g k, a_{***})$.
\medskip
Really, we have
\medskip
{\bf Conjecture 4.4.} All $M$ (uniformizable, having dual) satisfy Condition 4.2.
\medskip
{\bf Theorem 4.5.} Let $M$ satisfy Condition 4.2. Then $\g q'_M=\g q_{M'}$.
\medskip
{\bf Proof.} We denote the basis $\vf_*$ of $L(M')$ by $\hat \vf$ (matrix column). We have $\hat \eta= \g I\cdot \hat \vf$ where $\g I$ is a matrix of the change of bases. It is a skew $k'_*$-block anti-identity matrix $\left(\matrix 0&0&\dots&0&I_{k_{\g m+1}}\\ 0&0 & \dots&I_{k_{\g m}}&0\\ \dots & \dots & \dots & \dots & \dots \\ 0&I_{k_2}& \dots & 0&0 \\ I_{k_1}& 0& \dots & 0&0 \endmatrix \right)$, i.e. its antidiagonal block entries are identity matrices: the $(\g m+2-i,i)$-block is $I_{k_i}$, and other block entries are 0.
\medskip
We denote by $\Psi'_N$, resp. $\Psi'_{N,\eta}$ the $\th$-shift of the scattering matrix of $M'$ in bases $\hat \vf$, resp. $\hat \eta$ (as a base of $M'$ over $\p[T]$ we use $\hat f'$ in both cases). We have $\Psi'_N=\Psi_N^{t-1}\Xi_N^{-\g m}$ (Proposition 2.2), $\Psi'_{N,\eta}=\Psi'_N\g I^t$. Since $M$ satisfies Condition 4.2, there exists a set of Siegel matrices for $M'$ with respect to the basis $\hat \eta$. We denote it by $U_{****}$. It defines $B(U_{****})$ --- the corresponding $B$ in $M_r(\p)[N]$.
\medskip
We denote $\Psi_N\cdot B(S_{****})$, resp. $\Psi'_{N,\eta}\cdot B(U_{****})$ by $\g Z(M)$, resp. $\g Z(M^d)$. We have $\g Z(M)\in M_r(\p[[N]])$, $\g Z(M^d)\in M_r(\p[[N]])$ (Lemma 3.39), hence
$$\g Z(M)^t\cdot \g Z(M^d)=B(S_{****})^t\cdot \Psi_N^t\cdot  \Psi'_N \cdot \g I^t \cdot B(U_{****})=$$ $$B(S_{****})^t\cdot  \g I^t \cdot B(U_{****}) \cdot \Xi_N^{-\g m}\in M_r(\p[[N]])$$ and hence
$$  B(S_{****})^t\cdot  \g I^t \cdot B(U_{****})\cdot  \g I^t  \in \Xi_N^{\g m}M_r(\p[[N]])$$

We have  $\g I^t\cdot  B(U_{****})\cdot  \g I^t$ is of the form $X$ of Lemma 3.33. Further, $\Xi_N^{\g m}\in N^\g mM_r(\p[[N]])$, hence $\g I^t \cdot B(U_{****})\cdot \g I^t=\bar B(S_{****})$ (Lemma 3.33). This means that $U_{uvyz}=\bar S_{uvyz}$. The theorem follows from Lemma 3.36. $\square$
\medskip
\medskip
{\bf References}

\nopagebreak
\medskip
[A] Anderson Greg W. $t$-motives. Duke Math. J. Volume 53, Number 2 (1986), 457-502.
\medskip
[B] Brion, Michel. Lectures on the geometry of flag varieties. Topics in cohomological studies of algebraic varieties, 33–85,
Trends Math., Birkhäuser, Basel, 2005.
\medskip
[G] Goss, David Basic structures of function
field arithmetic.
\medskip
[GL07] A. Grishkov, D. Logachev,  Duality of Anderson t-motives. https://arxiv.org/pdf/0711.1928.pdf
\medskip
[GL17] A. Grishkov, D. Logachev, Lattice map for Anderson t-motives: first approach. J. of Number Theory. 2017, vol. 180, p. 373 -- 402. Electronic publication: https://doi.org/10.1016/j.jnt.2017.04.004 \ \ \  https://arxiv.org/pdf/1109.0679.pdf
\medskip
[GL20] A. Grishkov, D. Logachev, Introduction to Anderson t-motives: a survey. https://arxiv.org/pdf/2008.10657.pdf
\medskip
[GL21] A. Grishkov, D. Logachev, $h^1\ne h_1$ for Anderson t-motives. J. of Number Theory. 2021, vol. 225, p. 59 -- 89.
https://arxiv.org/pdf/1807.08675.pdf
\medskip
[HJ] Hartl U.; Juschka A.-K. Pink's theory of Hodge structures and the Hodge conjecture over function fields. "t-motives: Hodge structures, transcendence and other motivic aspects", Editors G. B\"ockle, D. Goss, U. Hartl, M. Papanikolas, European Mathematical Society Congress Reports 2020, and https://arxiv.org/pdf/1607.01412.pdf
\medskip
[P] Richard Pink, Hodge structures over function fields. Universit\"at Mannheim.
Preprint. September 17, 1997.
\medskip
[R16] Richarz, Timo. Affine Grassmannians and Geometric Satake Equivalences. 
Int. Math. Res. Not. IMRN 2016, no. 12, 3717 -- 3767
\medskip
\enddocument